\documentclass[11pt]{article}
 \usepackage{graphicx,amsmath,amsfonts,amssymb,color,mathrsfs}

\usepackage{epsfig}

 \newcommand{\QED}{\hfill \thicklines \framebox(6.6,6.6)[l]{}}

\numberwithin{equation}{section}

 \newtheorem{assumption}{Assumption}
 
 \newtheorem{theorem}{Theorem}[section]
 \newtheorem{lemma}{Lemma}[section]

 \newtheorem{example}{Example}[section]

 \newtheorem{definition}{Definition}[section]

\newcommand{\eqnb}{\begin{eqnarray*}}
\newcommand{\eqne}{\end{eqnarray*}}

\renewcommand\Re{\operatorname{Re}}

\def\R{{\mathbb R}}

\def\Re {{\rm Re}\,}
\def\E{{\mathbb E}}
\def\P{{\mathbb P}}

\def\C{{\mathbb C}}

\def\beqlb{\begin{eqnarray}}\def\eeqlb{\end{eqnarray}}
 \def\beqnn{\begin{eqnarray*}}\def\eeqnn{\end{eqnarray*}}

\title{\Large \bf Kernel Method --- An Analytic Approach \\ for Tail Asymptotics in Stationary Probabilities \\
of 2-Dimensional Queueing Systems}
\author{
Yiqiang Q. Zhao \\
School of Mathematics and Statistics \\
Carleton University, Ottawa, ON Canada K1S 5B6}
\date{February, 2020 \\ revised January 2021}

\begin{document}
 \maketitle

\begin{abstract}
In this paper, we provide a review on the kernel method, which is one of the options for characterizing so-called exact tail asymptotic properties in stationary probabilities of two-dimensional random walks, discrete or continuous (or mixed), in the quarter plane. Many two-dimensional queueing systems can be modelled via these types of random walks. Stationary probabilities are one of the most sought statistical quantities in queueing analysis. However, explicit expressions are available only for a very limited number of models. Therefore, tail asymptotic properties become more important, since they provide insightful information into the structure of the tail probabilities, and often lead to approximations, performance bounds, algorithms, among possible others.

Characterizing tail asymptotics for random walks in the quarter plane is a fundamental and also classical problem. Classical approaches are usually based on a complete determination of the transformation for the unknown probabilities of interest, for example, a singular integral presentation for the unknown probability generating function through boundary value problems \cite{FKM:82,Guillemin-Leeuwaarden:09}. In contrast to classical approaches (approaches based on the solution for the unknown probabilities or the transform of the unknow probabilities),
the kernel method, reviewed here, is very efficient for two-dimensional problems, which only requires the local information about the location of the dominant singularity of the unknown transformation function and the asymptotic property, through asymptotic analysis in complex analysis, at the dominant singularity.

This kernel method reviewed in this paper is an extension of the classical one, first introduced by Knuth~\cite{Knuth:69} and further developed, more than 30 years after, by Banderier \textit{et al.}~\cite{B-BM-D-F-G-GB:02}, targeting at one-dimensional problems. The method combines analytic continuation with asymptotic analysis (for example, see \cite{FIM:2017,Flajolet-Sedgewick:09}).
\vspace*{5mm}

\medskip


\end{abstract}

\section{Introduction}

The kernel method discussed here is an analytic method for tail asymptotics in a probability sequence, or distribution, which is an extended version developed from the classical kernel method. This method is a combination of analytic continuation and asymptotic analysis of complex functions, frequently used in analytic combinatorics. The idea of the classical kernel method can be traced back to Knuth~\cite{Knuth:69} in 1960s, which was further developed as the kernel method by Banderier \textit{et al.}~\cite{B-BM-D-F-G-GB:02} in early of this (the 21st) century.
To illustrate the key idea in the classical kernel method, suppose that we are interested in the unknown functions $F(x, y)$ and $G(x)$, which satisfy the following functional equation, often referred to as the fundamental form:
\begin{equation} \label{eqn:FM1}
    K(x, y)F(x, y) = A(x, y)G(x) + B(x, y),
\end{equation}
where $K(x,y)$, called the kernel function, $A(x,y)$ and $B(x,y)$ are given. The key idea, in the classical kernel method, for the solution of the unknown functions is quite intuitive. Consider pairs of $(x,y)$ satisfying the kernel equation $K(x,y)=0$. If $F(x,y)$ is finite, we then have
\begin{equation} \label{eqn:RHS-FM1-a}
    A(x, y)G(x) + B(x, y)=0.
\end{equation}
To solve the above equation for the unknown function $G(x)$, we use the fact that under some general conditions, $K(x,y)=0$ defines an algebraic (or multi-valued) function $y=Y(x)$ (similarly, $x=X(y)$). Select an analytic branch, say $Y_0(x)$, to rewrite the equation in (\ref{eqn:RHS-FM1-a}) as
\begin{equation} \label{eqn:RHS-FM1}
    A(x, Y_0(x))G(x) + B(x, Y_0(x))=0.
\end{equation}
This leads to the determination of
\[
    G(x)=-B(x, Y_0(x))/A(x, Y_0(x)),
\]
if $A(x, Y_0(x)) \neq 0$. Analytically continue $G(x)$ to a domain such that $(x, Y_0(x))$ is not a zero of the kernel function $K(x,y)$ and
substitute the continued version $G(x)$ into the fundamental form to have an expression or solution for $F$ given by
\[
    F(x, y) = \frac{-A(x, y)B(x,Y_0(x))/A(x, Y_0(x)) + B(x, y)}{K(x, y)}.
\]
Since the equation in (\ref{eqn:RHS-FM1}) contains only one unknown function, the associated problem is referred to as a problem with one unknown, or one-dimensional problem.

The focus of this paper is an extended version of the classical kernel method, which is combination of analytic continuation and asymptotic analysis of the unknown functions. This method is for studying stationary tail behaviour of two-dimensional queueing systems, or more generally for two-dimensional stochastic models. Assume that such a queueing system is modeled as a two-dimensional Markovian system. Under the stability condition, various stationary probability distributions (such as joint, marginal, boundary, directional among possible others) exist. We are interested in the transformations $\pi(x,y)$, $\pi_1(x)$ and $\pi_2(y)$ for an unknown joint probability sequence (or function), and two one-dimensional unknown probability sequences (or functions), respectively, which satisfy the following fundamental form:
\begin{equation} \label{eqn:FM2}
    -h(x,y)\pi(x,y)=h_{1}(x,y)\pi_1(x)+h_{2}(x,y)\pi_2(y)+h_{0}(x,y)\pi_{0,0},
\end{equation}
where $h(x,y)$, $h_1(x,y)$, $h_2(x,y)$ and $h_0(x,y)$ are given, and $\pi_{0,0}$ is an unknown constant. Applying the kernel method idea leads to
\begin{equation} \label{eqn:RHS-FM2-xy}
    h_{1}(x,y)\pi_1(x)+h_{2}(x,y)\pi_2(y)+h_{0}(x,y)\pi_{0,0}=0.
\end{equation}
Unlike the case with one unknown (or one-dimensional problem), the equation in (\ref{eqn:RHS-FM2-xy}) is a problem with two unknowns (functions) of two different variable $x$ and $y$, respectively. This type of problems is referred to as two-dimensional problem.
The classical kernel method leads to an equation in $x$:
\begin{equation} \label{eqn:RHS-FM2-x}
    h_{1}(x,Y_0(x))\pi_1(x)+h_{2}(x,Y_0(x))\pi_2(Y_0(x))+h_{0}(x,Y_0(x))\pi_{0,0}=0,
\end{equation}
or an equation in $y$:
\begin{equation} \label{eqn:RHS-FM2-y}
    h_{1}(X_0(y),y)\pi_1(X_0(y))+h_{2}(X_0(y),y) \pi_2(y)+h_{0}(X_0(y),y)\pi_{0,0}=0,
\end{equation}
but neither of them could immediately give a determination of $\pi_1(x)$ or $\pi_2(y)$.
In fact, solving for the unknowns $\pi_1(x)$ and $\pi_2(y)$ through the kernel equation, is very challenging, which has been very well documented in the literature (see our literature review section).
This process is equivalent to the dimension reduction from two to one, and the solution of resulting one dimensional problem.
Solving or determine $\pi_1(x)$ and $\pi_2(y)$ (therefore $\pi(x,y)$) is not the purpose of the kernel method reviewed here. The extended version of the kernel method focuses on exact tail asymptotics in stationary probabilities, or the stationary tail behaviour of the system, through studying their transformations, say $\pi_1(x)$, $\pi_2(y)$, $\pi(x,y)$, or possibly others, without having a complete characterization of the unknown functions. To achieve this goal, only the so-called local information about the unknown function, which is the location of the dominant singularity (or singularities) of the unknown function and its (or their) detailed asymptotic property as the variable approaches to the dominant singularity is required. This is actually a nature shared by other asymptotic methods, but not the classical methods in queueing literature. Very closely related to the kernel method reviewed here are asymptotic studies for queueing systems through a classical Tauberian theorem. In the kernel method, the connection of the asymptotic property of the unknown function to the tail asymptotic property of the probability function is established through the so-called Tauberian-like theorem (see, for example, Theorem~\ref{tauberian-1} and Theorem~\ref{5-thm1} in later sections). The key difference is that the monotone condition of the unknown probabilities required by the classical Tauberian theorem cannot be (or easily) verified, while the additional analyticity property of the transformation function for the unknown probabilities is often hold and can be verified for two-dimensional models.

There are two types of tail asymptotics, referred to as rough and exact in the literature. In this paper, we focus on exact tail asymptotics. Exact is used in contrast to rough, which is very often used when a large deviation decay is referred to. As a concrete example, suppose that the tail probability is equivalent to $n^{-1/2} \theta^n$ as $n \to \infty$ with $0<\theta<1$, which is called exact tail asymptotic, or exact tail decay rate, while the large deviation principle can only determine the rough decay rate $\theta$ (not the prefactor $n^{-1/2}$ in this example). Two nonnegative functions $f(t)$ and $g(t)$ are said to be tail equivalent, if $\lim_{t \to \infty} f(t)/g(t) =1$, denoted by $f(t) \sim g(t)$. In the discrete case, we are interested in probability sequences, say $f(k)$, and in the continuous case, we often take $f(t)$ to be a tail distribution, or $f(t)=\bar{F}(t)=1-F(t)$, where $F(t)$ is a distribution function. Both $f(k)$ and $F(t)$ can be non-proper: $\sum_k f(k) <1$ and $F(\infty)< 1$, respectively. For example, for the model described by the fundamental form in (\ref{eqn:FM2}), we are interested in the (exact) tail asymptotic of the probability function, whose transformation is $\pi_1(x)$, $\pi_2(y)$, $\pi(x,y)$, or other probabilities, say a marginal distribution of the joint distribution whose transformation is $\pi(x,y)$. We should point out that the tail asymptotic for the probability function whose transformation is $\pi_1(x)$ (or equivalently $\pi_2(y)$) is crucial for the tail asymptotic of other probability functions. Therefore, in this paper, we mainly concentrate on discussions of $\pi_1(x)$ (or $\pi_2(y)$).
 Due to symmetry, we can simply concentrate on one sequence (function). Specifically, let the probability sequence (function) be $f(t)$, we are aiming at finding a nonnegative function of form $g(t)= t^a \theta^t$ with $a \notin \{-1, -2, \ldots \}$ and $0<\theta<1$ such that
\[
    f(t) \sim g(t).
\]
There are two main reasons for us to consider the form $t^a \theta^t$ of tail asymptotics: in the analysis we need to use Tauberian-like theorems that require this type of asymptotics, and for two-dimensional queueing systems we have not found any other types of tail asymptotics reported in the literature.

There are three key components in the kernel method: (1) The solution of the kernel equation. This can be achieved through various approaches, geometric, analytic or algebraic. We will treat the solution as an algebraic (two-valued) function and characterize these two branches in terms of complex analysis at the basic level to avoid high level mathematical concepts. The analytic property of these two functions are crucial for the kernel method; (2) Dimension reduction to one-dimensional problems. This is the main additional challenge when extending the classical kernel method to the extended version reviewed in this paper. This is achieved through the interplay of the two branches determined in (1) and the interplay of the two unknown transformation functions; and (3) Asymptotic analysis at the dominant singularity. Since the dimension has been reduced from two to one in (2), the problem becomes classical involving characterization of the detailed asymptotic behaviour at the dominant singularity of the unknown (one variable) function, and applications of the Tauberian-like theorem.

The purpose of this paper is to review the key components in the kernel method and to show how the kernel method
can be used for exact tail asymptotic properties of various two-dimensional or related models.
To this end, we divides our discussions in four steps for each of the cases, consisting of:
\begin{description}
\item[Step~1:] Fundamental form: formulated from the system dynamics;
\item[Step~2:] Analytic continuation of the unknown function: through various available methods;
\item[Spet~3:] Asymptotic analysis: identifying the location of the dominant singularity, and detailing the asymptotic property at the dominant singularity based on the interplay of the two unknown functions, say $\pi_1(x)$ and $\pi_2(y)$ defined in (\ref{eqn:RHS-FM2-x}) and (\ref{eqn:RHS-FM2-y}));
\item[Step~4:] A proper Tauberian-like theorem: connecting the asymptotic property of the unknown function to the tail asymptotic property of the probability function whose transformation is the unknown function.
\end{description}

The rest of this paper consists of five sections. In the next section, Section~\ref{sec:2}, we provide a literature review. In Sections~\ref{sec:3}, \ref{sec:4} and \ref{sec:5}, we discuss the main asymptotic results for random walks in the quarter plane (the discrete case), for a Brownian model (the continuous case) and for a fluid model modulated by an infinite-state Markov chain (a mixed case), respectively, and provide details for the four steps. In Section~\ref{sec:6}, we discuss possible extensions of the kernel method.

\section{Alternative methods and literature review} \label{sec:2}

Two-dimensional queueing models play an important role in applications. For example, models, such as join-the-shortest-queue systems, arrivals with two separate demands, priority queues, modified Jackson networks, find many applications in various fields. It is well known that except for some special cases (say, models with product-form solutions), in general, we do not expect to have explicit formulas for the joint stationary probability distribution for the model. Instead, approximation methods are often a tool for studying such a model. Two-dimensional queueing systems having discrete states are often modeled as a random walk in the quarter plane. Continuous state queueing systems are also important, either as an approximate model for queueing systems, or as a limiting model for queueing systems (say diffusion approximations), which can be modeled as continuous-state random walks in the quarter plane.

For two-dimensional problems, approaches focusing on characterizing the transform of an unknown probability sequence or function include boundary-value problems, analytic or algebraic methods, uniformization approach, among possible others. Literature in this area is vast, including:
Kingman~\cite{Kingman:1961}, who is among the earliest researchers using a functional equation of two variables for studying two-dimensional queueing systems. In \cite{Kingman:1961}, the symmetric join-the-shortest-queue model was considered. Complex analysis, analytic continuation, series expansion are the main tools for his study. Tail asymptotic properties were obtained based on the determination of the transformation functions of probabilities of interest, which is referred to as a classical method in this review. 
Malyshev~\cite{Malyshev:72,Malyshev:73}, whose work started a systematic treatment of two-dimensional queueing systems in term of analytic approach. In \cite{Malyshev:72}, a fundamental form for the unknown generating functions of the (general) random walk in the quarter plane was obtained, which is essentially the same as (\ref{eqn:FM2}). Analytic continuation of the unknown generating functions was discussed through studying the kernel equation. His work was then further developed by Fayolle and Iasnogorodski into the book \cite{FIM:2017}, which contains systematic studies for the random walk in the quarter plane, through algebraic methods and boundary value problems. Many of our Section~\ref{sec:3} results are from this book. Tail asymptotic properties were obtained for the so-called simple random walk (random walks in the quarter plane without transitions to north-east, north-west, south-east or south-west). For this special type of random walk, the author showed that all branch points are real, which is the base of the same claim for the general random walk in the quarter plane (see \cite{FIM:2017}). The tail asymptotic property was obtained through a saddle point approximation method on the integral expression of the unknown generating function, which is also classified as a classical method. Other related work to \cite{FIM:2017} include 
    Fayolle and Iasnogorodski~\cite{FI:1979}, 
    Fayolle, King and Mitrani~\cite{FKM:82}, 
    Kurkova and Suhov~\cite{Kurkova-Suhov:03},
    Bousquet-Melou~\cite{Bousquet-Melou:05},
     Raschel~\cite{Raschel:10},
     Kurkova and Raschel~\cite{Kurkova-Raschel:2011,Kurkova-Raschel:2013},
which are only a small sample of the large number of related references. 
Flatto and McKean~\cite{Flatto-McKean:77} and Flatto and Hahn~\cite{Flatto-Hahn:84} used a different method, referred to the uniformization method for explicit solutions of two-dimensional queueing systems. The same symmetric join-the-shortest-queue model considered in \cite{Kingman:1961} was considered in \cite{Flatto-McKean:77}, which is a walk of genus 0, not discussed in this review. The two-demand queueing model (which is used as an example in Section~\ref{sec:3} of this paper) was considered in \cite{Flatto-Hahn:84}. This is a random walk of genus 1. Analytic continuation and the explicit solution of the unknown generating functions are based on the property of elliptic functions. This is the first example for which the explicit solution is not rational, but algebraic. The characterization of the solution to the kernel equation is equivalent to the characterization through universal covering (see \cite{FIM:2017}) of the Riemann surfaces defined by the kernel equation, but only basic analysis concepts are involved in the uniformization method. Therefore, people in the queueing area, who are not very familiar with concepts of Reimann surfaces, could prefer this method.
This two-demand model was extended to the generalized two-demand model, which was first considered by 
Wright~\cite{Wright:92}, using a similar approach. Integral presentations for the unknown transformations were obtained, based on which some tail asymptotic results were presented but with errors (see, \cite{Li-Zhao:10b}). 
For boundary value problems, the monograph by Cohen and Boxma~\cite{cohen-boxma:83} is a good reference for queueing models. Exact tail asymptotic analysis based on boundary value problems are also considered classical, since it is based on the solution of the unknown transformation functions, Related work includes 
Guillemin, Knessl and Leeuwaarden~\cite{GKL2011},
Guillemin and Pinchon~\cite{Guillemin-Pinchon:2004} and Guillemin and Leeuwaarden~\cite{Guillemin-Leeuwaarden:09}.

Our focus in this paper is on exact tail asympotics in terms of the kernel method.
Other approaches (not discussed in detail in this paper) for tail asymptotics of two-dimensional models are available, such as compensation method, methods based on Markov additive processes, large deviations, among possible others.
For example, the compensation method was proposed by Adan~\cite{AWZ:1993}, which provides a complete solution of the joint probability distribution in terms of an infinite series of product-form terms.
Studies based on properties of Markov additive processes (including matrix-analytic methods) are vast, including
    McDonald~\cite{McDonald:99},
    Takahashi, Fujimoto and Makimoto~\cite{TFM:01},,
    Foley and McDonald~\cite{Foley-McDonald:01,Foley-McDonald:05a,Foley-McDonald:05b} (the exact tail asymptotics are based on detailed counting of the green function),
    Haque~\cite{Haque:03},
    Kroese, Scheinhardt and Taylor~\cite{KST:04},    
    Miyazawa~\cite{Miyazawa:04},    
    Miyazawa and Zhao~\cite{Miyazawa-Zhao:04},
    Haque, Liu and Zhao~\cite{Haque-Liu-Zhao:05},
    Li and Zhao~\cite{Li-Zhao:05},
    Motyer and Taylor~\cite{Motyer-Taylor:06},     
    Li, Miyazawa and Zhao~\cite{Li-Miyazawa-Zhao:07},
    He, Li and Zhao~\cite{He-Li-Zhao:08},
    Liu, Miyazawa and Zhao~\cite{Liu-Miyazawa-Zhao:08}, 
    Tang and Zhao~\cite{Tang-Zhao:08}, 
    Adan, Foley and McDonald~\cite{Adan-Foley-McDonald:09},  
    Khanchi~\cite{Khanchi:08,Khanchi:09},
        Kobayashi, Miyazawa and Zhao~\cite{Kobayashi-Miyazawa-Zhao:10},
        Kobayashi and Miyazawa~\cite{Kobayashi-Miyazawa:2011}.
Not just for the rough decay, this type of method can also lead to conditions for exact geometric decay. However, for other types of exact tail asymptotics, detailed counting on the Green's function, or asymptotic expansions of transformation functions were involved.
Examples of using large deviations approach can be found in  Borovkov and Mogul'skii~\cite{Borovkov-M:01} and Lieshout and Mandjes~\cite{LM:2008},.
    A method based on geometric properties of the model, initiated by Miyazawa, is robust including            
        Miyazawa~\cite{Miyazawa:07,Miyazawa:09,Miyazawa:08,Miyazawa:11},
        Miyazawa and Rolski~\cite{Miyazawa-Rolski:09},
        Ozawa~\cite{Ozawa:2013},
        Ozawa and Kobayashi~\cite{Ozawa-Kobayashi:2018}. 
When for exact tail asymptotics, the asymptotic analysis and the Tauberian-like theorem components were integrated into this geometric method.


The kernel method, reviewed in this paper, is an extension based on the work by Knuth~\cite{Knuth:69} and Banderier \textit{et al.}~\cite{B-BM-D-F-G-GB:02}. The Tauberian-like theorems are developed from the work in analytic combinatorics,  for which references can be found in Flajolet and Sedgewick~\cite{Flajolet-Sedgewick:09}. The extended version of the kernel method has been successfully applied to various models such as:
    Li and Zhao~\cite{Li-Zhao:10b} for the generalized two-demand queueing model;
    Li, Tavakoli and Zhao~\cite{Li-Tavakoli-Zhao:11} for genus 0 random walks in the quarter plane; 
    Ye~\cite{Ye:2012} and Ye, Li and Zhao~\cite{Ye-Li-Zhao:2015} for a longer-queue-serve-first system; 
    Zafari~\cite{Zafari:2012} for the generalized join-the-shortest-queue model considered earlier by several other researchers including \cite{Foley-McDonald:01}, \cite{Houtum-etc:2001}, \cite{Kurkova-Suhov:03} and \cite{Miyazawa:09};     
    Dai and Zhao~\cite{Dai-Zhao:2013} for a revisit of the wireless 3-hop networks with stealing considered in \cite{GKL2011};
    Song, Liu and Dai~\cite{Song-Liu-Dai:2015} for a discrete-time preemptive priority queueing model; 
    Song, Liu and Zhao~\cite{Song-Liu-Zhao:2015} for a revisit of the retrial queue with two input streams and two orbits considered by Avrachenkov, Nain and Yechiali~\cite{ANY:2014};
    Dai, Dawson and Zhao~\cite{Dai-Dawson-Zhao:2015} for extending the kernel method to continuous random walks;
    Dai, Kong and Song~\cite{Dai-Kong-Song:2017} for a two-stage queue;
    Dai, Dawson and Zhao~\cite{Dai-Dawson-Zhao:2018} for extending the kernel method to a 3-dimensional tandem queueing system; 
    Li and Zhao~\cite{Li-Zhao:2018} for a systematic treatment of the kernel method for random walks in the quarter plane;
    Li, Liu and Zhao~\cite{Li-Liu-Zhao:2019} for a 2-demand model modulated by a two-state Markov chain, which is a special case of the random walk in the quarter plane modulate by a finite-state Markov chain;
    Song and Lu~\cite{Song-Lu:2021} for the Israeli queue with retrials and non-persistent customers.

It should be emphasized that the above discussed references are far away from a complete list of related work.

\section{Two-dimensional discrete-state systems} \label{sec:3}

The kernel method, discussed in this paper, for two-dimensional queueing systems started from discrete-state systems. The earliest report was provided by Li and Zhao~\cite{Li-Zhao:10b}, on the generalized 2-demand model. The same model was considered in Wright~\cite{Wright:92}, for which an integral presentation of the unknown generating function was obtained. Since \cite{Li-Zhao:10b}, the kernel method has been used for exact tail asymptotics of other queueing models (see the previous section for more details).
A systematic study of the kernel method on the tail asymptotic for random walks in the quarter plane is available in Li and Zhao~\cite{Li-Zhao:2018}.
This is the same type of model was studied in the book by Fayolle, Iasnogorodski and Malyshev~\cite{FIM:2017}, in which more references can be found. The book contains algebraic methods and boundary value problems for characterizing the whole unknown generating functions. The kernel method discussed in this paper is closely related to their work. In fact, the analytic continuation of the unknown generating functions was a direct result from their book.

In this section, we provide a summary on the exact tail asymptotic results obtained in terms of the kernel method with discussions, remarks and comments. We also provide an example to demonstrate the key steps of using the kernel method.

The random walk in the quarter plane is a discrete-time Markov chain with state space $S=\{ (m,n) | m, n =0, 1, 2, \ldots \}$ and its (nonzero) transition probabilities $P_{(m,n), (m',n')}$ are described by:
\begin{equation}
    P_{(m,n), (m+i,n+j)} = \left \{ \begin{array}{ll}
p^{(0)}_{i,j}, & \text{if } (m,n)=(0,0), \\
p^{(1)}_{i,j}, & \text{if } (m,n)=(m,0) \text{ with } m \geq 1, \\
p^{(2)}_{i,j}, & \text{if } (m,n)=(0,n) \text{ with } n \geq 1, \\
p_{i,j}, & \text{if } m \geq 1 \text{ and } n \geq 1,
\end{array} \right.
\end{equation}
where $p^{(0)}_{i,j}$, $p^{(1)}_{i,j}$, $p^{(2)}_{i,j}$ and
$p_{i,j}$ are arbitrary non-negative real numbers satisfying
\[
    \sum_{i, j =0, 1} p^{(0)}_{i,j} =1, \;
    \sum_{i=0, \pm 1, j =0, 1} p^{(1)}_{i,j} =1, \;
    \sum_{i =0, 1, y=\pm 1} p^{(2)}_{i,j} =1, \;
    \sum_{i, j =0, \pm 1} p^{(0)}_{i,j} =1.
\]
This is a two-dimensional random walk model with reflective horizontal and vertical  boundaries. At state $(m,n)$, the chain could only move to one of its eight neighboring states, or remain in $(m,n)$ in the next transition, constrained to the quarter plane. Therefore, all possible values for $i$ and $j$ are defined correspondingly.
Assume that the Markov chain is irreducible, aperiodic and positive recurrent, under which $\pi_{m,n}$ is the unique stationary joint probability vector for the model.
In this paper, we only focus on the light-tailed behaviour, which is equivalent to assume that $M =(M_{x,} M_{y})= ( \sum_i i ( \sum_j p_{i,j} ), \sum_j j  ( \sum_i p_{i,j}  ) ) \neq 0$ (see, for example, Lemma~3.3 of \cite{Kobayashi-Miyazawa:2011}).

We now discuss the four key steps in applying the kernel method and provide the main asymptotic results with discussions. The is assumed stable. For a stable condition, one may refer to Theorem~3.3.1 of \cite{FIM:2017}, which has been amended by Lemma~2.1 in \cite{Kobayashi-Miyazawa:2011}.

\textbf{Step~1:} To state the fundamental form, we follow the book \cite{FIM:2017} to define the following generating functions for the probability sequences associated with the interior states, horizontal boundary states and vertical boundary states, respectively:
\begin{align}
    \pi(x,y) &= \sum_{m=1}^{\infty} \sum_{n=1}^{\infty} \pi_{m,n}x^{m-1}y^{n-1}, \\
    \pi_1(x) &=\sum_{m=1}^{\infty} \pi_{m,0}x^{m-1}, \\
    \pi_2(y) &=\sum_{n=1}^{\infty} \pi_{0,n}y^{n-1}.
\end{align}
The probability generating functions (PGF) $\pi(x,y)$, $\pi_1(x)$ and $\pi_2(y)$ are respectively referred to as the PGF of the joint probabilities $\pi_{m,n}$ excluding the boundaries, the PGF of the horizontal boundary probabilities excluding the origin, and the PGF of the vertical  boundary probabilities excluding the origin.
According to equation (1.3.6) in~\cite{FIM:2017}, the fundamental form for the random walk in the quarter plane is given by (\ref{eqn:FM2}), that is,
\begin{equation} \label{eqn:FM}
    -h(x,y)\pi(x,y)=h_{1}(x,y)\pi_1(x)+h_{2}(x,y)\pi_2(y)+h_{0}(x,y)\pi_{0,0},
\end{equation}
where
\begin{align}
    h(x,y) &= xy\left( \sum_{i=-1}^{1}\sum_{j=-1}^{1}p_{i,j}x^{i}y^{j}-1\right), \\
    h_{1}(x,y) &=x\left(\sum_{i=-1}^{1}\sum_{j=0}^{1}p_{i,j}^{(1)}x^{i}y^{j}-1\right), \\
    h_{2}(x,y) &=y\left(\sum_{i=0}^{1}\sum_{j=-1}^{1}p_{i,j}^{(2)}x^{i}y^{j}-1\right), \\
    h_{0}(x,y) &=\left(\sum_{i=0}^{1}\sum_{j=0}^{1}p_{i,j}^{(0)}x^{i}y^{j}-1\right).
\end{align}
There are a few optional ways for obtaining the fundamental form. For a specific model, such a functional equation can be obtained directly by elementary calculations based on the system of equations for the model, which involves collecting coefficients of $z^n$. This direct method could also work for the general model, but tedious calculations and arrangements are involved, which is not an efficient approach. The probabilistic argument employed in \cite{FIM:2017} is obviously a preferable approach, which is much neater. Their method can be used for obtaining the fundamental form for other types of models, for example for random walks in the quarter plane modulated by a finite-state Markov chain
(see Section~\ref{sec:6}).
\bigskip

\textbf{Step~2:} There are a few options for the analytic continuation of the unknown function $\pi_1(x)$ (and $\pi_2(y)$ by symmetry), including direct methods, the uniformization method, the method through the universal covering among others. Once the analytic continuation of the functions $\pi_1(x)$ and $\pi_2(y)$ is established, the analytic continuation for $\pi(x,y)$ can be easily discussed.

In this step, it is crucial to study the kernel equation: $h(x,y)=0$.  The random walk in the quarter plane is called non-singular if $h(x,y)$, as a polynomial in the two variables $x$ and $y$ over real numbers, is irreducible and quadratic in both variables. We will focus on the non-singular random walks only in this paper, since the analysis for a singular random walk is either easier or can be done similarly to that for the non-singular walk (see, for example, Li, Tavakoli and Zhao~\cite{Li-Tavakoli-Zhao:11}). For the non-singular random walk, given $x$, we write $h(x,y)$ as a quadratic form in $y$:
 \[
    h(x,y)= a(x)y^{2}+b(x)y+c(x),
 \]
where
\begin{align*}
a(x) =&p_{-1,1}+p_{0,1}x+p_{1,1}x^{2}, \\
b(x) =&p_{-1,0}-(1-p_{0,0})x+p_{1,0}x^{2}, \\
c(x) =&p_{-1,-1}+p_{0,-1}x+p_{1,-1}x^{2}.
\end{align*}
 Then, there are two solutions for $y$ according to the quadratic formula, given by
 \[
       Y_{\pm}(x) = \frac{-b(x) \pm \sqrt{D_1(x)}}{2a(x)},
 \]
where $D_{1}(x)=b^{2}(x)-4a(x)c(x)$ is the discriminant. For a given $y$, we can have a similar discussion, and the discriminant is denoted by $D_2(y)$.
So, the kernel equation defines an algebraic (or two-valued) function $Y(x)$.
If $D_1(x)=0$ does not have multiple roots, the random walk is called a genus one walk (in this case, $D_2(y)=0$ does not have multiple roots either), otherwise, the random walk is called a genus zero walk. The genus one walk represents majority of queueing systems and it is more difficult to analyze. Therefore, in this paper, we focus on the genus one walk only.

It is well-known that the complex square root function $\sqrt{z}$, which is the branch such that it is equal to the real square root function when $z \geq 0$ is real, is not analytic on the whole complex $z$-plane $\mathbb{C}_{z}$, but it is analytic on the cut plane $\widetilde{\mathbb{C}}_{z} =\mathbb{C}_{z} \setminus (-\infty,0]$, where $z=0$ is the branch point of the square root function. For the case of random walks, each root $x$ of $D_1(x)=0$ is called a branch point of $Y(x)$. For the non-singular genus one random walk,
we know that all four branch points $x_i$ ($i=1, 2, 3, 4$) are real, satisfying:
\begin{equation} \label{eqn:branch-points}
    |x_1| < x_2 < 1 < x_3 < |x_4|,
\end{equation}
if $M_y \neq 0$ (for example, see Lemma~2.3.8 of \cite{FIM:2017}). It will be helpful if we could find an interpretation or intuitive explanation about why all branch points are real. It is believed that this should be a property of the random variable, but such an interpretation is still sought for. It is known that for a function expressed as a Taylor series of nonnegative coefficients, the radius of convergence should be a dominant singularity of the function according to Pringsheim's Theorem. We then immediately know there are at most two complex branch points. However, no simple ways available for the result in (\ref{eqn:branch-points}). Explicit criteria given in terms of $p_{m,n}$ on when $x_1$ and $x_4$ are positive are also given in Lemma~2.3.8 of \cite{FIM:2017}, but they are not crucial in our discussions in this paper.

The (two-valued) function $Y(x)$ plays a crucial role in the kernel method. The analytic behaviou of $Y(x)$ is then a key property for our analysis. Since given a value of $x$, $Y(x)$ equals to $Y_+(x)$ or $Y_-(x)$, we should first exam the property of these two functions, or the square root function $\sqrt{D_1(x)}$. Similar to the square root function $\sqrt{z}$, we consider the following cut plane: $\widetilde{\widetilde{\mathbb{C}}}_{x} =\mathbb{C}_{x} \setminus [x_{3},x_{4}] \cup [x_{1},x_{2}]$, where $[x_{3},x_{4}] = [x_3, \infty) \cup (-\infty, x_4]$ if $x_4 <0$. Unfortunately, neither $Y_-(x)$ nor $Y_+(x)$ is analytic in this cut plane.   Both of them are analytic on domains of the cut plan, separated by curves defined by $D_1(x) =0$. For a specific example, one may refer to the 2-dimensional tandem queue model studied in \cite{Guillemin-Leeuwaarden:09}. For this tandem queue model, the authors in \cite{Guillemin-Leeuwaarden:09} demonstrated how to construct, using the building blocks $Y_+$ and $Y_-$, an analytic continued function on the whole cut plan. For a general case considered in this paper, we can construct the analytic continued function $Y_0(x)$, which equals to either $Y_-(x)$ or $Y_+(x)$, according to $|Y_0(x)| = \min( |Y_-(x)|,|Y_+(x)|)$ on the cut plane $\widetilde{\widetilde{\mathbb{C}}}_{x}$ (for example, see \cite{Li-Zhao:2018}). In the book  \cite{FIM:2017}, $Y_0$ was constructed using a different way (more advanced mathematics involved) and this function was proved to be $Y_-$ or $Y_+$ whichever has the minimum modulus. According to the uniqueness property of the analytic continued function, our $Y_0$ coincides with the $Y_0$ in \cite{FIM:2017}.
The functions $Y_0(x)$ and $X_0(y)$ (by symmetry), and the interplays defined by (\ref{eqn:RHS-FM2-x}) and (\ref{eqn:RHS-FM2-y}) play an important role in analytic continuation of the unknown function $\pi_1(x)$ (and $\pi_2(y)$).

A few options are available for the analytic continuation of $\pi_1(x)$. For example, a similar direct argument to that in \cite{Guillemin-Leeuwaarden:09} was used for a queueing system with two coupled processors in Fayolle and Iasnogorodski~\cite{FI:1979}. A uniformization method was used by Flatto and Hahn~\cite{Flatto-Hahn:84} for the 2-demand model in terms of properties of elliptic functions.
In terms of the universal covering and Riemann surface theory, a nicer proof for a general random walk in the quarter plane was presented in \cite{FIM:2017}. It is worth mentioning that the uniformization method is essentially equivalent to the universal covering approach, the former is mathematically less demanding while the latter requires more advanced mathematical concepts. The analytic continuation results is presented in the following lemma.

\begin{lemma}[Theorem~3.2.3 in \cite{FIM:2017}] \label{lemma1.3}
For the stable non-singular, genus one random walk, $\pi_1(x)$ is a meromorphic function in the complex cut plane $\widetilde{\mathbb{C}}_{x}=\mathbb{C}_{x} \setminus [x_{3},x_{4}]$. Similarly, $\pi_2(y)$ is a meromorphic function in the complex cut plane $\widetilde{\mathbb{C}}_{y}=\mathbb{C}_{y} \setminus [y_{3},y_{4}]$.
\end{lemma}
\bigskip

\textbf{Step~3:} There are three keys things involved in this step: (a) the location of the dominant singularity (or singularities); (b) the type of a dominant singularity; and (c) the detailed asymptotic property at a dominant singularity. Besides many other aspects, the interplay of the two unknown functions plays a key role in the analysis. The interplay was briefly mentioned in the introduction (see equations (\ref{eqn:RHS-FM2-x}) and (\ref{eqn:RHS-FM2-y})), and to be precise, we state the following property.

\begin{theorem}[Theorem~4.3 in  \cite{Li-Zhao:2018}] \label{theorem1.1}
In the cut plane $\widetilde{\widetilde{\mathbb{C}}}_{x}$,
\begin{equation} \label{eqn:1.9}
    \pi_1(x)=\frac{-h_{2}(x,Y_{0}(x))\pi_2(Y_{0}(x))-h_{0}(x,Y_{0}(x))\pi_{0,0}}{h_{1}(x,Y_{0}(x))},
\end{equation}
except at a zero of $h_{1}(x,Y_{0}(x))$, or at a pole of $\pi_1(x)$ or $\pi_2(Y_{0}(x))$.

Similarly, in the cut plane $\widetilde{\widetilde{\mathbb{C}}}_{y}$,
\begin{equation} \label{eqn:1.11}
    \pi_2(y)=\frac{-h_{1}(X_{0}(y),y)\pi_1(X_{0}(y))-h_{0}(X_{0}(y),y)\pi_{0,0}}{h_{2}(X_{0}(y),y)},
\end{equation}
except at a zero of $h_{2}(X_{0}(y),y)$, or at a pole of $\pi_2(y)$ or $\pi_1(X_{0}(y))$.
\end{theorem}
This interplay is important for deriving properties, including the above theorem, of $\pi_1(x)$ as a function of $x$ (and $\pi_2(y)$ as a function of $y$), or equivalent to dimension reduction from two $(x,y)$ to one $x$ or $y$. The dimension reduction is achieved by introducing the uniformized variable in the uniformization, and through the generator $\delta = \eta \xi$ of two Galois automorphisms $\eta$ and $\xi$ in the algebraic method in \cite{FIM:2017}.

From the above theorem, it is clear that the radius of convergence of $\pi_1(x)$ is either $x_3 >1$ or $1 \leq |x_p| \leq x_3$, where $x_3$ is a branch point at which $\pi_1(x)$ is not analytic, and $x_p$ is a pole of $\pi_1(x)$, since for the generating function of a nonnegative sequence, the radius $x_p$ of convergence should be a dominant singularity, $\pi_1(x)$ is meromorphic in the cut plane, and the model is light-tailed. To make our discussion simpler, we exclude the X-shaped random walks. Specifically, we assume that $p_{i,j}$ is not X-shaped, where $p_{i,j}$ is called X-shaped if $p_{i,j}=0$ for all $i$ and $j$ such that $|i+j|=1$ (see, for example, \cite{Li-Zhao:2018}).
It was proved in \cite{Li-Zhao:2018} that for any non-X-shaped random walk, there exists a unique dominant singularity on the convergence radius, while for the X-shaped random walk, there are exactly two dominant singularities. All tail asymptotic properties for the non-X-shaped random walks hold for the X-shaped random walk, but in a periodic manner.

Based on the above discussions, the dominant singularity of $\pi_1(x)$ is either the branch point $x_3$ or the unique minimum pole $1< x_p \leq x_3$.
The following theorems characterize the nature of the dominant singularity of $\pi_1(x)$.

\begin{theorem}[Theorems~4.4, 4.5 and 4.7 in  \cite{Li-Zhao:2018}]  \label{theorem1.2}
 Let $x_p$ be a pole of $\pi_1(x)$ with the
smallest modulus greater than one. Assume that $|x_p| \leq x_3$. Then, one of the
following two cases must hold:

\textbf{1.} $x_p$ is a zero of $h_{1}(x,Y_{0}(x))$. In this case, $x_p \stackrel{\triangle}{=}x^*$ is the unique zero of $h_{1}(x,Y_{0}(x))$ in $(1,x_3]$;

\textbf{2.} $\widetilde{y}_0=Y_{0}(x_p)$ is a zero of
$h_{2}(X_{0}(y),y)$ and $|\widetilde{y}_0|>1$. In this case, $x_p\stackrel{\triangle}{=}\widetilde{x}_1=X_1(\widetilde{y}_0)$, where $\widetilde{y}_0$ is the unique
zero of the function $h_{2}(X_{0}(y),y)$ in $(1,y_{3}]$.
\end{theorem}

Now, it is clear that the location of the dominant singularity is at $R$, the radius of convergence of $\pi_1(x)$ satisfying
\[
    R = x_{dom} = \min (x_3, x^*, \widetilde{x}_1).
\]
The type of asymptotic property at the dominant singularity depends on which one (or which ones since the branch point and the pole can be equal) of the three candidates equals (or equal) the convergence radius $R$. For convenience, we let $x^*=\infty$ if there does not exist a pole $x^* \leq x_3$. A similar convention is made for $\widetilde{x}_1$. From a standard asymptotic analysis, we can detail the asymptotic property of $\pi_1(x)$ for all possible cases. Such a detailed analysis confirms that
a total of four cases exist regarding the type of the dominant singularity.
\begin{description}
\item[Case~1.] $x_{dom}=\min\{x^{\ast}, \widetilde{x}_{1}\}<x_{3}$
with $x^{\ast} \neq \widetilde{x}_{1}$, or
$x_{dom}=\widetilde{x}_{1}=x^{\ast}=x_{3}$;

\item[Case~2.] $x_{dom}=x_{3}=\min\{x^{\ast},\widetilde{x}_{1}\}$ with $x^{\ast}\neq
\widetilde{x}_{1}$;

\item[Case~3.] $x_{dom}=x_{3}<\min\{x^{\ast},\widetilde{x}_{1}\}$;

\item[Case~4.] $x_{dom}=x^{\ast}=\widetilde{x}_{1}<x_{3}$.
\end{description}

Corresponding to these four cases, the four types of asymptotic properties at the dominant singularity are detailed in the following theorem.
\begin{theorem}[Theorems~4.8 in  \cite{Li-Zhao:2018}] \label{theorem3.1}
For the function $\pi_1(x)$, a total of four types of
asymptotic properties exist as $x$ approaches to a dominant singularity of
$\pi_1(x)$, corresponding to the type of the dominant singularity.

\textbf{Case~1.}
\begin{equation*}
    \lim_{x\rightarrow x_{dom}}\left( 1-\frac{x}{x_{dom}}\right) \pi_1(x)=c_{0,1}(x_{dom});
\end{equation*}

\textbf{Case~2.}
\begin{equation*}
    \lim_{x \to x_{dom}}\sqrt{1-x/x_{dom}}\pi_1(x)=c_{0,2}(x_{dom});
\end{equation*}

\textbf{Case~3.}
\begin{equation*}
    \lim_{x\rightarrow x_{dom}}\sqrt{1-x/x_{dom}}\pi_1^{\prime}(x)=c_{0,3}(x_{dom});
\end{equation*}

\textbf{Case~4.}
\begin{equation*}
    \lim_{x\rightarrow x_{dom}}\left( 1-\frac{x}{x_{dom}}\right)
    ^{2}\pi_1(x)=c_{0,4}(x_{dom}),
\end{equation*}
where $c_{0,i}(x_{dom})$, $i=1, 2, 3, 4$, are non-zero constant.
\end{theorem}
It should be point out that the non-zero constant result is needed when we apply the Tauberian-like theorem. For a specific model, the proof could be much easier by using the specific model structure, which the proof for the general case is much more difficult.

There always exists an explicit expression for $x_3$ in terms of the system parameters $\{ p_{i,j}, p^{(k)}_{i,j}; k = 0, 1, 2 \}$ since it is a root of a polynomial of degree 4. However, in most cases, explicit expressions for the pole may not be available. The main reason is the structure of the expression for $Y_0$ or $X_0$, which is rather complex. To solve equation $h_1(x, Y_0(x)) = 0$ for $x$, one commonly used method is through polynomialization of $h_1$, which is always feasible. The resulting polynomial is often of a degree higher than 4. Depending on models, for some of which we could find a solution explicitly, but for most of which numerical procedures need to be developed for computations. This (for roots of polynomials) can be efficiently achieved. The two-demand model reviewed in this paper is one of the examples in the literature, which has a nice explicit expression for the pole.
\bigskip

\textbf{Step~4:} The final step of the kernel method is the application of a Tauberian-like theorem, which bridges the asymptotic property of the unknown function $\pi_1(x)$ at its dominant singularity and the tail asymptotic property of the coefficients, or $\pi_{m,0}$, of $\pi_1(x)$. For the discrete-state case, like the random walks in the quarter plane, the following Tauberian-like theorem is an immediate consequence of Corollary~VI.1 in \cite{Flajolet-Sedgewick:09}.

\begin{definition}[Definition~VI.1 in Flajolet and Sedgewick~\cite{Flajolet-Sedgewick:09}]
For given numbers $\varepsilon >0$ and $\phi$ with $0<\phi <\pi/2$, the open domain
$\Delta (\phi, \varepsilon)$ is defined by
\begin{equation}
    \Delta (\phi, \varepsilon) = \left \{z \in \mathbb{C}: |z| < 1+\varepsilon, z \neq 1, \arg |z-1| > \phi \right \}.
\end{equation}
A domain is a $\Delta$-domain at 1 if it is a $\Delta(\phi, \varepsilon)$ for some $\varepsilon>0$ and $0< \phi< \pi/2$. For a complex number $\zeta \neq 0$,
a $\Delta$-domain at $\zeta$ is defined as the image $\zeta \cdot \Delta(\phi, \varepsilon)$ of a $\Delta$-domain $\Delta(\phi, \varepsilon)$ at 1
under the mapping $z \mapsto \zeta z$. A function is called $\Delta$-analytic if it is
analytic in some $\Delta$-domain.
\end{definition}

\begin{figure}[h]
\centering
      \includegraphics[width=8cm]{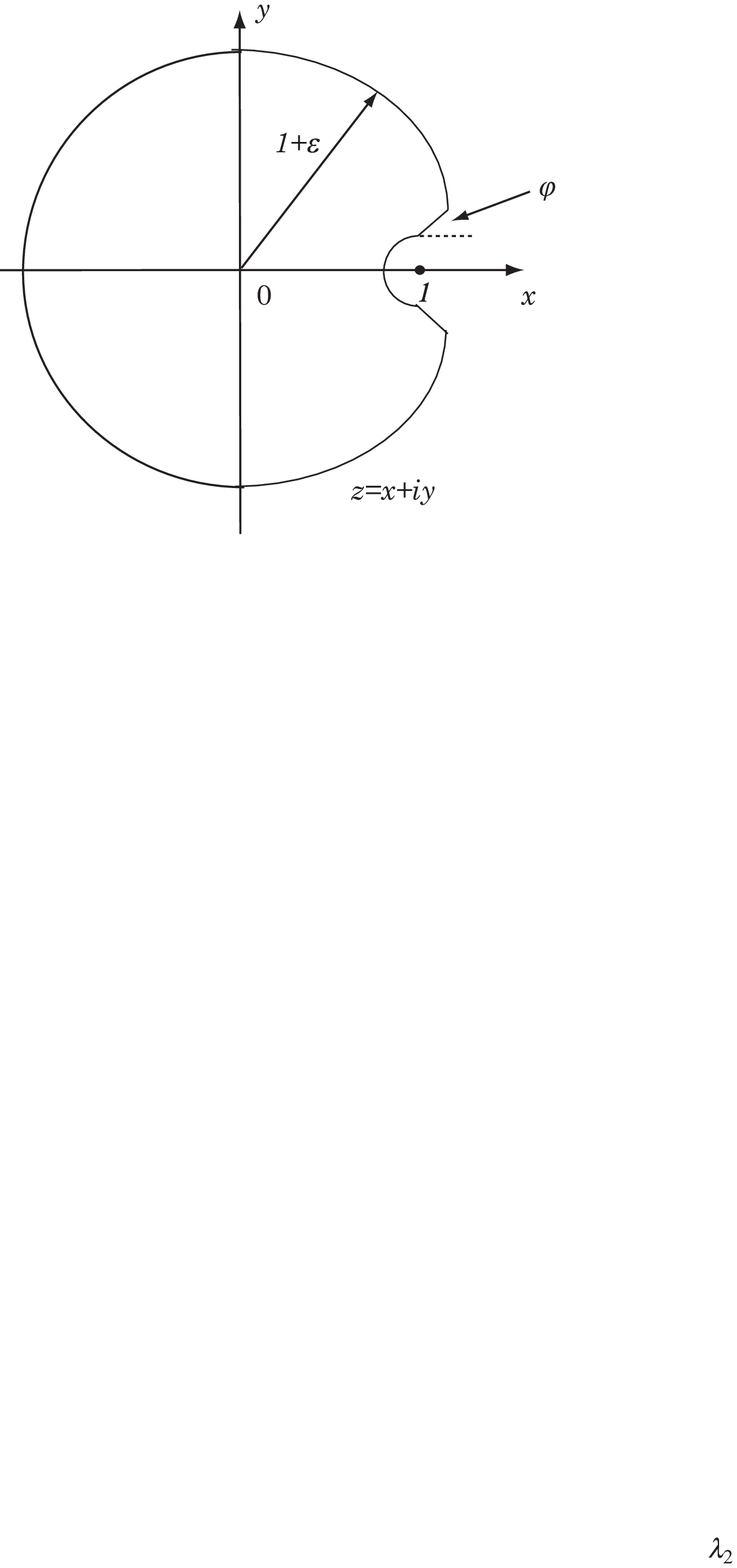}
\vspace*{-11cm}
     \caption{$\Delta$-domain.}
\end{figure}

\begin{theorem}[Tauberian-like theorem for single singularity] \label{tauberian-1}
Let $A(z)=\sum_{n\geq 0}a_{n}z^{n}$ be analytic at 0 with $R$ being the
radius of convergence. Suppose that $R$ is a singularity of $A(z)$
on the circle of convergence such that $A(z)$ can be continued to
a $\Delta$-domain at $R$. If for a real number $\alpha \notin \{0,
-1, -2, \ldots\}$,
\begin{equation} \label{eqn:constant}
    \lim_{z \rightarrow R}(1-z/R)^{\alpha}A(z)= g,
\end{equation}
where $g$ is a non-zero constant, then,
\begin{equation*}
    a_{n} \sim  \frac{g}{\Gamma(\alpha)}  n^{\alpha-1} R^{-n},
\end{equation*}
where $\Gamma(\alpha)$ is the value of the gamma function at $\alpha$.
\end{theorem}

Figure~1 provides a depiction of the $\Delta$-domain. This domain is required by the Tauberian-like theorem. It is called Tauberian-like to distinguish it from the classical Tauberian theorem. Both the Tauberian-like theorem and its corresponding classical Tauberian theorem (for example, see Feller~\cite{Feller:1971}) provide the same asymptotic conclusion, but with different required conditions. $a_n$ being eventually monotone ($a_n$ is monotone for $n > N$ where $N$ is finite) is required by classical Tauberian theorems, which cannot be verified (or easily verified) since $a_n$ is the unknown probability sequence sought for. This is the main drawback for applying the classical Tauberian theorem, including using Heaviside operational calculus. The effort of removing this condition from the classical Tauberian theorem started in analytic combinatorics, which led to the kernel method reviewed in this paper. As a tradeoff, more on analyticity of the function $A(z)$ is required by the Tauberian-like theorem, which requires that $A(z)$ can be analytically continued to a $\Delta$-domain. This can be verified for the random walk in the quarter plan and other models discussed in this paper. Applying this extended version of the kernel method to queueing models for tail asymptotic properties was initiated in \cite{Li-Zhao:10b}, and then used by other researchers.
The above Tauberian-like theorem is a version for a single dominant singularity on the radius of convergence circle, but this is not crucial. Tauberian-like theorems for multiple dominant singularities can be found in \cite{Flajolet-Sedgewick:09}.

In order to see that $\pi_1$ is analytic in a $\Delta$-domain, we first notice Lemma~\ref{lemma1.3}, which is Theorem~3.2.3 in \cite{FIM:2017} and poles can only be discrete isolated. Since we only consider the light-tailed case (under the condition of $M \neq 0$) in this review, the radius $R$ of convergence of $\pi_1(x)$ is greater than one.
For non-singular genus 1, non-X-shaped random walks, there is only one, which is $R$, dominant singularity (see \cite{Li-Zhao:2018} for details). Therefore, $\pi_1(x)$ ($A(x)$ in the Tauberian-like theorem) can be analytically continued to the $\Delta$-domain at $R>1$.

The application of the Tauberian-like theorem to our case is straightforward,  which leads to the main tail asymptotic result for the two-dimensional random walks.
It is worth mentioning here that in Case~3, the asymptotic property presented in Theorem~\ref{theorem3.1} is for the derivative of $\pi_1$. When applying the Tauberian-like theorem, we notice that
\[
    \lim_{z \rightarrow 1} \frac{f(z)}{K(1-z)^s} = \lim_{z \rightarrow 1} \frac{f'(z)}{K s(1-z)^{s-1}}
\]
according to L'H\^{o}pital's rule.
All involved limits in the following theorem, and also in Theorem~\ref{theorem3.1}, are taken over the $\Delta$-domain.
\begin{theorem}[Theorems~5.1 in  \cite{Li-Zhao:2018}]
\label{theorem4.1-a} Consider the stable non-singular, non-X-shaped genus 1
random walk. Corresponding to the above four cases, we have the
following tail asymptotic properties for the boundary
probabilities $\pi_{n,0}$ for large $n$. In all cases,
$c_{0,i}(x_{dom})$ ($1\leq i\leq 4$) are given in
Theorem~\ref{theorem3.1}.

\textbf{Case~1.} (Exact geometric decay)
\begin{equation} \label{eqn:exact}
    \pi_{n,0} \sim c_{0,1}(x_{dom})\left( \frac{1}{x_{dom}}\right)^{n-1};
\end{equation}

\textbf{Case~2.} (Geometric decay multiplied by a factor of
$n^{-1/2}$)
\begin{equation}
    \pi_{n,0} \sim \frac{c_{0,2}(x_{dom})}{\sqrt{\pi}}n^{-1/2}\left(\frac{1}{x_{dom}}\right)^{n-1};  \label{eqn:1/2}
\end{equation}

\textbf{Case~3.} (Geometric decay multiplied by a factor of
$n^{-3/2}$)
\begin{equation} \label{eqn:3/2}
    \pi_{n,0} \sim \frac{c_{0,3}(x_{dom})}{\sqrt{\pi}}n^{-3/2}\left(\frac{1}{x_{dom}}\right)^{n-1};
\end{equation}

\textbf{Case~4.} (Geometric decay multiplied by a factor of $n$)
\begin{equation}
    \pi_{n,0} \sim c_{0,4}(x_{dom})n\left( \frac{1}{x_{dom}}\right)^{n-1}. \label{eqn:n}
\end{equation}
\end{theorem}

By symmetry, the tail asymptotic properties of $\pi_{0,n}$ can be obtained without any extra effort. Based on the asymptotic properties for $\pi_{n,0}$ and $\pi_{0,n}$, we can further discuss tail asymptotics for $\pi_{n,m}$ and $\pi_{m,n}$ for any fixed $m$ and $n$, respectively, for the marginal distributions, and for the probability distribution along a diagonal (non-coordinate) direction (see, for example, \cite{Li-Zhao:2018}).

By concluding the discussion for exact tail asymptotic properties for the random walk in the quarter plane, we provide the 2-demand model as an example.

\begin{example}[\cite{Flatto-Hahn:84}] \rm
This 2-demand mode was studied first by Flatto and Hahn in \cite{Flatto-Hahn:84}. This is a specific case of the random walk in the quarter plane, which is non-singular and genus 1. This parallel system consists of two independent exponential servers with service rates $\mu_1$ and $\mu_2$, respectively. In front of each of the two servers,  there is a queue of infinite capacity. The Poisson arrival process with rate $\lambda$ to the system simultaneously creates two jobs, one to each queue. Jobs in each queue
are processed according to the first-come-first-served discipline.
Arrivals and services are assumed to be independent. If we use the number of jobs in each queue (including the job in the service) as the system variable, we then have a two-dimensional continuous-time Markov chain with states $(m,n)$, where $m$ and $n$ are the numbers of jobs in queue 1 and queue 2, respectively.
A necessary and sufficient stability condition for this system is $\lambda < \mu_i$ for $i=1, 2$. Without loss of generality, we use $\lambda + \mu_1 +\mu_2 =1$ for the uniformization parameter to convert this continuous-time Markov chain to a discrete-time Markov chain with the same stationary probability vector $\pi_{m,n}$.

For the random walk in the quarter plane, sometime it is more convenient to define the following slightly different generating functions:
\[
    P(x,y) = \sum_{m, n=0}^{\infty} \pi_{m,n} x^m y^n,
\]
$P_1(x)=P(x,0)$ and $P_2(y)=P(0,y)$.
It is clear that $P(x,y)$ is the PGF of the complete joint distribution $\pi_{m,n}$ (including the probabilities on the boundaries), and $P_1(x)$ and $P_2(y)$ are the PGFs for the boundary probabilities, horizontal and vertical respectively, including the probability $\pi_{0,0}$ at the origin.
We then have the following slight modified fundamental form:
\begin{equation} \label{eqn:fundamental}
    H(x,y)P(x,y)=H_{1}(x,y)P_1(x)+H_{2}(x,y)P_2(y)+ H_{0}(x,y) \pi_{0,0},
\end{equation}%
where
\begin{align*}
    H(x,y) &= -h(x,y), \\
    H_1(x,y) &= -h(x,y)+h_1(x,y)y, \\
    H_2(x,y) &= -h(x,y) + h_2(x,y) x, \\
    H_0(x,y) &=  h_0(x,y) x y + h(x,y) - h_1(x,y)y - h_2(x,y)x.
\end{align*}
(\ref{eqn:fundamental}) can be easily obtained based on the fundamental form in (\ref{eqn:FM}).
For the 2-demand model, we specifically have
\begin{align*}
H(x,y) &=-xy+ \mu_{1}y + \mu_{2}x + \lambda x^{2}y^{2}, \\
H_{1}(x,y) &= \mu_{2}x(1-y), \\
H_{2}(x,y) &=\mu_{1}y(1-x), \\
H_{0}(x,y) &=0.
\end{align*}
Here, $H_0 =0$ simplifies the analysis significantly.
For this model, an explicit expression for the unknown function $P_1(x)$ was obtained in terms of a uniformization method, based on which exact tail asymptotic properties for $\pi_{n,0}$ were reported in \cite{Flatto-Hahn:84}. Here, without a determination of the unknown function, we use its local information about the location and type of the dominant singularity for the same tail asymptotic results.

\begin{theorem}[Theorem~4.3 in {Li-Zhao:10b}] \label{the:3}
For the 2-demand model, let $\hat{x}=\mu_1/\lambda$. We can show that:

\textbf{(i)} $\hat{x} \leq x_3$, and the equality holds if and only if $\mu_1=\mu_2=\mu$;

\textbf{(ii)} $\hat{x}$ is the unique zero of $H_{1}(x,Y_{0}(x))$ in $(1,x_3]$ if and only if $\mu_1 \leq \mu_{2}$;

\textbf{(iii)} $H_{2}(X_{0}(Y_{0}(x)),Y_{0}(x))$ has no solution in $(1,x_{3}] $ with $|Y_{0}(x)|>1$.
\end{theorem}
Therefore, we only have two candidates, $x_3$ and $x^* = \hat{x}$, for the dominant singularity $x_{dom}$, for which the following three cases exist:
\begin{description}
\item[Case~1.] $\mu_1 < \mu_2$. In this case, the dominant singularity $x_d$ is a pole equal to $x_{dom}= x^* = \hat{x} < x_3$;

\item[Case~2.] $\mu_1 = \mu_2 =\mu$. In this case, $x_{dom}=x_{3}=x^* = \hat{x}= \mu/\lambda$. Both the pole and the branch point are the dominant singularity;

\item[Case~3.] $\mu_1 > \mu_2$. In this case, no pole exits in $(1,x_3]$ and $x_{3}$ is the dominant singularity.
\end{description}
The exact tail asymptotic is characterized by

\begin{theorem} \label{the:4}
For the boundary probabilities $\pi_{m,0}$ of the 2-demand queueing model, the exact tail asymptotic is given by:
\begin{description}
\item[Case~1.]
\[
    \pi_{m,0} \sim c_{1} \left( \frac{\lambda}{\mu_{1}} \right)^{m};
\]

\item[Case~2.]
\[
    \pi _{m,0} \sim c_{2} \; m^{-1/2} \left(\frac{\lambda}{\mu_{1}}\right)^{m};
\]

\item[Case~3.]
\[
    \pi_{m,0} \sim c_{3} \; m^{-3/2} \left(\frac{1}{x_{3}}\right)^{m},
\]
where, $P_2(1) = 1- (\lambda/\mu_1)$,
\begin{align*}
    c_{1}&=\frac{(\mu_{2}-\lambda \hat{x}) P_2(1)}{\mu_{2}}, \\
    c_2 &= \frac{\mu_{1}(\hat{x}-1) P_2(1)}{\sqrt{\mu_{1}\mu_{2}(\hat{x}-x_{1})(\hat{x}-x_{2})}\sqrt{\pi}}, \\
    c_{3} &= \frac{(x_{3}-1)\mu_{1} \sqrt{4\mu_{2}\lambda x_{3}(x_{3}-x_{1})(x_{3}-x_{2})}}{4 \mu_{2}\lambda x_{3}^{3}\sqrt{\pi}}
 \left[ \frac{P_2(Y_{0}(x_{3})) + Y_{0}(x_{3})(1-Y_{0}(x_{3})) P_2^{\prime}(Y_{0}(x_{3}))}{(Y_{0}(x_{3})-1)^{2}}\right].
\end{align*}
\end{description}
\end{theorem}

\end{example}

Exact tail asymptotic characterization for a generalized 2-demand queueing model can be found in \cite{Li-Zhao:10b}.

\section{Two-dimensional continuous-state systems} \label{sec:4}

The two-dimensional semimartingale reflective Brownian motions (SRBM), as the counterpart to the random walk in the quarter plane, is a type of continuous-state Markovian system (random walk), which plays a fundamental role in both theoretical and applied issues (see, for example, Dai and Harrison~\cite{DH1992} and Williams~\cite{W1995,W1996}), including in queueing applications.
In this section, we explain how the key ideas in the kernel method for the discrete random walk can be extended for continuous random walks, based on the work in Dai, Dawson and Zhao~\cite{Dai-Dawson-Zhao:2015}. Since the key components of the kernel method applied to continuous models remain the same as applied to discrete models, our discussions will focus more on the differences.   

The SRBM is a process $Z(t)=(Z_1(t),Z_2(t))$ defined by $Z(0)=X(0) \in \mathbb{R}^2_+ =\{(x_1,x_2); x_1, x_2 \geq 0\}$ and for $t>0$,
 \begin{equation} \label{def1}
    Z(t)=X(t)+R Y(t),
 \end{equation}
where $X(t)$ is a Brownian motion (random walk) with the drift vector $\mu$ and the covariance matrix $\Sigma=(\Sigma_{ij})$,  $R=(r_{ij})$ is $2 \times 2$ matrix, and
$Y(t)=(Y_1(t), Y_2(t))$ is a $2$-dimensional process such that
\begin{enumerate}
\item[(i)] $Y(t)$ is continuous and non-decreasing with $Y(0)=0$;

\item[(ii)] $Y_j(t)$ only increases at times $t$ for which $Z_j(t)=0$, $j=1, 2$;

\item[(iii)] $Z(t)\in\mathbb{R}^2_+$, $t\geq 0$.
\end{enumerate}
It follows from the above definition that the SRBM $Z(t)$ behaves like an ordinary Brownian motion, and reflects, when it hits a boundary, to the direction specified by the corresponding column of $R$ and regulated by the process $Y(t)$.

We assume that the process is stable, or $R$ is non-singular and $R^{-1}\mu<0$:
\[
    r_{11}>0,\;r_{22}>0,\;\textrm{and}\;r_{11}r_{22}-r_{12}r_{21}>0;
\]
and
\[
    r_{22}\mu_1-r_{12}\mu_2<0,\;\textrm{and}\; r_{11}\mu_2-r_{21}\mu_1<0,
\]
which is a necessary and sufficient condition. Under this condition, let $Z=(Z_1,Z_2)$ be a random vector that has
the stationary distribution $\pi$ of the SRBM, let $Z(0)$ follow the stationary distribution $\pi$, and let $\E_\pi(\cdot)$
denote the conditional expectation given that $Z(0)$ follows the stationary distribution $\pi$.

The tail property of the joint stationary distribution $\pi$ is the focus of the analysis. For this purpose, we consider the transformation of $\pi$ and set up a proper (functional) equation, the fundamental form, which links the joint distribution $\pi$ (through its transformation) to a couple of one-dimensional (simpler) distributions. In the discrete random walk case, the boundary probabilities are natural choices. However, for the continuous random walk case, this is no longer valid (since the probability on a line is simple zero). We need to define proper boundary measures such that through the system dynamics, such that a link (fundamental form) can be established. A standard approach is to define a measure on the quarter plane by the expected occupation time according to a local time, which is the regulation process $Y(t)$ in our case. 
Following Dai and Harrison~\cite{DH1992}, we define
 \begin{equation} \label{M-2}
    V_i(A)=\E_\pi\Big[\int_{0}^11_{\{Z(u)\in A\}}dY_i(u)\Big], \quad i=1,2,
 \end{equation}
where $A\subset\mathbb{R}^2_+$ is a Borel set. According to (\ref{M-2}), $V_i$ defines a finite measure on $\mathbb{R}_+^2$, since $\E_\pi(Y(1))$ is finite, and has a support on the face $F_i=\{x=(x_1,x_2)\in \mathbb{R}^2_+:x_i=0\}$.
This boundary measure is the natural counterpart to the boundary probability vectors in the discrete case. Now, we can define the transformations for the joint probability measure and the two boundary measures:
\begin{align*}
    &\phi(\theta_1,\theta_2) =\mathbb{E}_\pi e^{\langle \theta,Z\rangle}, \\
    &\phi_1(\theta_2) =\int_{\mathbb{R}^2_+}e^{\theta_2x_2}V_1(dx)=\E_{\pi}\int_{0}^1e^{\theta_2Z_2(u)}dY_1(u),\\
&\phi_2(\theta_1) =
\int_{\mathbb{R}^2_+}e^{\theta_1x_1}V_2(dx)=\E_{\pi}\int_{0}^1e^{\theta_1Z_1(u)}dY_2(u).
\end{align*}
\bigskip

\textbf{Step~1:}
The fundamental form follows directly from (2.3) in Dai and Miyazawa~\cite{DM2011}:
\begin{equation}\label{2-6}
    \gamma(x,y)\phi(x,y)=\gamma_1(x,y)\phi_1(y)+\gamma_2(x,y)\phi_2(x),
\end{equation}
where
\begin{align}\label{2-7}
    \gamma_1(x,y) & =r_{11}x+r_{21}y, \\
\label{2-5}
    \gamma_2(x,y) &=r_{12}x+r_{22}y, \\
\label{2-4}
    \gamma(x,y) &=-<\hat{x},\mu>-\frac{1}{2}<\hat{x},\Sigma \hat{x}> \\
    &=x\mu_1+y\mu_2+\frac{1}{2}\Sigma_{11}x^2+\Sigma_{12}xy+\frac{1}{2}\Sigma_{22}y^2\nonumber
\\&=\frac{1}{2}\Sigma_{22}y^2+(\mu_2+\Sigma_{12}x)y+\frac{1}{2}\Sigma_{11}x^2+x\mu_1\nonumber
\\&= a y^2+b(x)y+c(x).
\end{align}
The proof of (\ref{2-6}) was based on a stand basic adjoint relationship in SRBM (for example, Harrison and Williams~\cite{Harrison-Williams:1987} and Dai and Harrison~\cite{DH1992}). 

It is worthwhile to mention that the boundary equations given by $\gamma_i(x,y)=$ for $i=1, 2$ define two straight lines (when both $x$ and $y$ are real), which are simpler than that for the discrete case, which are quadratic curves.
\bigskip

\textbf{Step~2:} For analytic continuation of $\phi_2(x)$, we consider the kernel equation $\gamma(x,y)=a y^2+b(x)y+c(x)=0$, which is a quadratic form in $y$ for fixed $x$ with the two solutions given by
\begin{equation} \label{A-9}
    Y_{\pm}(x)=\frac{-b(x)\pm\sqrt{D_1(x)}}{2a},
\end{equation}
where $D_1(x)=b^2(x)-4ac(x)$ is the discriminant. The branch points play an important role in the continuation.

\begin{lemma}[Lemma~1 in \cite{Dai-Dawson-Zhao:2015}] \label{3-lem1}
$D_1(x)$  has two real zeros $x_1$ and $x_2$ satisfying $x_1\leq0<x_2$. Furthermore,
$D_1(x)>0$ in $(x_1,x_2)$, and $D_1(x)<0$ in $(-\infty,
x_1)\cup(x_2,\infty)$.
\end{lemma}
Based on the above property and some direct arguments, we can show (see, for example, \cite{Dai-Dawson-Zhao:2015}) that:
\begin{lemma}[Lemma~2 in \cite{Dai-Dawson-Zhao:2015}]\label{2-lem2}
Both   $Y_{+}(x)$ and $Y_{-}(x)$  are analytic on the cut plane $\widetilde{\mathbb{C}}_x= \mathbb{C}_x \setminus\big\{(-\infty,x_1]\cup[x_2,\infty)\}$.
\end{lemma}

We recall that for the discrete case discussed in the previous section, both $Y_+(x)$ and $Y_-(x)$ cannot be analytic in the whole cut plane. For the continuous case, it is much simpler to construct the analytic branches $Y_0(s)$ (and $X_0(y)$). It can be verified that $Y_0(x) = Y_-(x)$ (and $X_0(y)=X_-(y)$). We use the analytic branch $Y_0(x)$ in the continuation of the function $\phi_2(x)$. Once again, the interplay between the two unknown functions $\phi_1(y)$ and $\phi_2(x)$ is the key for the continuation.
Instead of the radius of convergence, we define the convergence parameter for the moment generating function here. For a probability distribution $F(x)$ on $\R_+$, let $g(x)=\int_{0}^\infty e^{\lambda x}dF(x)$ be the moment generating function, whose convergence parameter $C_p(g)$ is defined by
\[
   C_p(g)=\sup\{\lambda\geq 0: g(\lambda)<\infty\}.
\]
We now allow $x$ to be complex. Then, the Pringsheim's theorem (see, for example, Dai and Miyazawa~\cite{DM2011} and Markushevich~\cite{M1977}) says that $g(z)$ is analytic on $\{z\in \mathbb{C}_z: \Re (z)<C_p(g)\}$.

Let $\tau_1=C_p(\phi_2)$. We can then prove in terms of the interplay (see Lemma~6 in  \cite{Dai-Dawson-Zhao:2015}) that $\tau_1 >0$. It implies that the boundary measure $V_2$ is light-tailed. Through some additional direct arguments, we can show that:
\begin{lemma}[Lemmas~7 and 10 in \cite{Dai-Dawson-Zhao:2015}] \label{4-lem4}
$\phi_2$ can be analytically continued to the region: $\{x \in \widetilde{\mathbb{C}}_x: \gamma_2(x,Y_0(x))\neq 0\} \cap \{x \in \widetilde{\mathbb{C}}_x: \Re (Y_0(x))<\tau_2\}$, and
\beqlb\label{4-2}
\phi_2(x)=-\frac{\gamma_1\big(x,Y_0(x)\big)\phi_1\big(Y_0(x)\big)}{\gamma_2\big(x,Y_0(x)\big)},
\eeqlb
where $\tau_2=C_p(\phi_1)$.
Furthermore, the function $\phi_2(x)$ is  meromorphic on the cut plane $\widetilde{\mathbb{C}}_x$.
\end{lemma}

\bigskip

\textbf{Step~3:} Following similar arguments to the discrete case, we have the following property about the pole:
\begin{theorem}[Lemmas~9 and 11 in \cite{Dai-Dawson-Zhao:2015}]  
 Let $x_p$ be a pole of $\phi_2(x)$ with the largest modulus smaller than or equal to $x_2$. Assume that $|x_p| >0$. Then, one of the
follow two cases must hold:

\textbf{1.} $x_p$ is a real zero of $\gamma_2(x,Y_{0}(x))$. This is true if and only if $\gamma_2(x,Y_{0}(x)) \geq 0$. In this case, $x_p =x^*$ is the unique zero of $\gamma_2(x,Y_{0}(x))$ in $(0,x_2]$;

\textbf{2.} $\widetilde{y}=Y_{0}(x_p)$ is a real zero of
$\gamma_1(X_{0}(y),y)$ and $\widetilde{y}>0$. In this case, $x_p=\widetilde{x}=X_1(\widetilde{y})$, where $\widetilde{y}$ is the unique
zero of the function $\gamma_1(X_{0}(y),y)$ in $(0,y_2]$.
\end{theorem}
Therefore, the same as for the discrete case, $x^*$ (a pole),  $\widetilde{x}$ (a pole), and $x_2$ (a branch point) are the three candidates for the dominant singularity of $\phi_2(x)$. 
We noticed that for the discrete case, $D_1(x)$ is a polynomial of degree 4, but for the continuous case, it is a polynomial of degree 2. It implies that the explicit expression for the branch point $x_2$ can be much simpler, obtained directly from the quadratic formula. In addition, an explicit expression for $x^*$ always exists (see equation (60) in \cite{Li-Zhao:2018}). For the general case, we may not have an explicit expression for $\widetilde{x}$.

Similar to the discussion for the discrete case, we have the following four cases for the dominant singularity of $\gamma_2(x)$:
\begin{description}
\item[Case~1.] $\tau_1=\min\{x^{\ast}, \widetilde{x}\}<x_{2}$
with $x^{\ast} \neq \widetilde{x}$, or
$\tau_1=\widetilde{x}=x^{\ast}=x_{2}$;

\item[Case~2.] $\tau_1=x_{2}=\min\{x^{\ast},\widetilde{x}\}$ with $x^{\ast}\neq
\widetilde{x}$;

\item[Case~3.] $\tau_1=x_{2}<\min\{x^{\ast},\widetilde{x}\}$;

\item[Case~4.] $\tau_1=x^{\ast}=\widetilde{x}<x_{2}$.
\end{description}

Without any difference, the detailed asymptotic property at the dominant singularity $\tau_1$ can be obtained through standard asymptotic analysis, which is described by the following theorem:
\begin{theorem}[Theorem~1 in \cite{Dai-Dawson-Zhao:2015}]  \label{thm1}
For the function $\phi_2(x)$, a total of four types of asymptotics exist as $x$ approaches to $\tau_1$.
\begin{description}
\item[Case~1.]
    \beqlb\label{thm1-1}
    \lim_{x\to \tau_1}(\tau_1-x)\phi_2(x)=A_1(\tau_1);
    \eeqlb

\item[Case~2.] \beqlb\label{thm1-2}
\lim_{x\to\tau_1}\sqrt{\tau_1-x}\phi_2(x)=A_2(\tau_1); \eeqlb

\item[Case~3.]
\beqlb\label{thm1-3}
\lim_{x\to\tau_1}\sqrt{\tau_1-x}\phi_2'(x)=A_3(\tau_1); \eeqlb

\item[Case~4.]
\beqlb\label{thm1-4}
\lim_{x\to\tau_1}(\tau_1-x)^2\phi_2(x)=A_4(\tau_1). \eeqlb
\end{description}
\end{theorem}
In the theorem, all constants $A_i(\tau_1)$ for $i=1,2,3,4$ can be explicitly expressed as shown in \cite{Dai-Dawson-Zhao:2015}.
\bigskip

\textbf{Step~4:} The SRBM discussed here was also considered in \cite{DM2012}. To obtain exact tail asymptotic results they combined the geometric method with the asymptotic analysis and a Tauberian-like theorem for the continuous case established in \cite{DM2012}. The following theorem is a slightly different version of the Tauberian-like theorem for the continuous case with weaker conditions.

Denote
 \beqnn
    \Delta(z_0,\epsilon)=\big\{z\in\C:z\neq z_0,\;|{\rm
arg}(z-z_0)|>\epsilon\big\},
 \eeqnn
 where ${\rm arg}(z) \in(-\pi,\;\pi]$ is the principal part of the argument of a complex number $z$.
\begin{theorem}[Theorem~2 in \cite{Dai-Dawson-Zhao:2015}] \label{5-thm1}
Assume that $g(z)=\int_{0}^\infty e^{st}f(t)dt$ satisfies the following conditions:
\begin{description}
\item[(1)]The left-most singularity of $g(z)$ is $\alpha_0$ with
$\alpha_0>0$. Furthermore, we assume that as $z\to \alpha_0$,
\[
    g(z)\sim(\alpha_0-z)^{-\lambda}
\]
for some $\lambda\in\C\setminus \{0, -1, -2, \ldots \}$;

\item[(2)]$g(z)$ is
analytic on $\Delta(\alpha_0,\epsilon_0)$ for some
$\epsilon_0\in(0,\frac{\pi}{2}]$;

\item[(3)] $g(z)$ is bounded on
$\Delta(\alpha_0,\epsilon_1)$ for some $\epsilon_1>0$.
\end{description}
Then, as $t\to\infty$,
\beqlb\label{R1-16}
f(t)\sim e^{-\alpha_0 t}\frac{t^{\lambda-1}}{\Gamma(\lambda)},
\eeqlb
where $\Gamma(\cdot)$ is the Gamma function.
\end{theorem}
A direct application of the above Tauberian-like theorem leads to the following tail asymptotic properties for the boundary measure $V_2$.

\begin{theorem}[Theorem~3 in \cite{Dai-Dawson-Zhao:2015}]
For the  boundary measure $V_2\big(x,\;\infty\big)$, we have the following tail asymptotic properties for large $x$:
\begin{description}
\item[Case~1.]
\[
    V_2\big(x,\infty\big)\sim C_1e^{-\tau_1 x};
\]
\item[Case~2.]
\[
    V_2\big(x,\infty\big)\sim C_2 x^{-\frac{1}{2}} e^{-\tau_1 x};
\]

\item[Case~3.]
\[
    V_2\big(x,\infty\big)\sim C_3 x^{-\frac{3}{2}} e^{-\tau_1 x};
\]

\item[Case~4.]
\[
    V_2\big(x,\infty\big)\sim C_4 x e^{-\tau_1 x},
\]
\end{description}
where $C_i$, $i=1,2,3,4$ are constants.
\end{theorem}

Before ending our discussions on the continuous model, we would like to emphasize that it is now well understood that for rough decay (and possible for exact geometric decay), various methods are available. However, for exact tail asymptotics, more detailed calculations should be involved, for example asymptotic expansions based on asymptotic analysis used in the kernel method, which was also incorporated into the geometric method (\cite{DM2012,Kobayashi-Miyazawa:2013}), and detailed counting for the green function (\cite{Foley-McDonald:01, Foley-McDonald:05b}).

\section{Two-dimensional systems with mixed discrete and continuous states} \label{sec:5}

In this section, we demonstrated, by a fluid model driven by the $M/M/c$ queue, how the kernel method can be extended to a case, for which one dimension is discrete (the number of customers in the $M/M/c$) and the other dimension is continuous (the fluid content). The content of this section is based on Li, Liu and Zhao~\cite{Li-Liu-Zhao:2019}.
Our discussions in this section will emphasize the differences and special challenges for the mixed type of models.

Let $X(t)$ be the fluid level at time $t$ in the fluid queueing system, driven by the $M/M/c$ queueing system $Z(t)$, the state or the number of the customers in the $M/M/c$ queueing system at time $t$. It means that the rate of change of the fluid level depends on the state $Z(t)$. Specifically, let $r_{Z(t)}$ denote the rate of change of the fluid level, or the net input rate, at time $t$. Then, the dynamics of the fluid level $X(t)$ is described by the following differential system:
\[
    \frac{dX(t)}{dt}=\left\{ \begin{array}{ll}
 r_{Z(t)}, & \mbox{if $X(t)> 0$ or $r_{Z(t)}\geq 0$,} \\
     0, & \mbox{if $X(t)= 0$   and $r_{Z(t)}< 0$,}
   \end{array}
   \right.
\]
where $r_{Z(t)}=r> 0$ if $Z(t) \geq c$, and $r_{Z(t)}=Z(t)-c < 0$ if $Z(t) \leq c-1$. The credit depletes the fluid during the partial busy period of an $M/M/c$ queue (i.e. whenever an arriving customer finds less than $c$ customers in the queue) at a negative rate $r_{Z(t)}$. It is reasonable to assume that the negative rate $r_{i}$ increases in $i$. Without loss of generality, we assume that  the net input rate is $r_{i}= i-c$ for any $0\leq i\leq c-1$.

We assume the $M/M/c$ queue $Z(t)$ with arrival rate $\lambda$ and service rate $\mu$ for each server is stable, or $\lambda< c\mu$. Let $\xi_{i}$, $i=0, 1, 2, \ldots$, be the unique stationary distribution for $Z(t)$ and let $Z$ be the stationary variable for $Z(t)$ (or $Z$ has the distribution $\xi_i$). For the fluid queue to be stable, we assume that $\sum_{i\in E}\xi_{i}r_{i}<0$. Let $X$ be the stationary variable for $X(t)$. We are interested in the exact tail asymptotic of the buffer content probability distribution
\[
    \Pi(x) = \sum_i \Pi_{i}(x),
\]
through the analysis of the joint distribution
\[
  \Pi_{i}(x)=P\{Z=i, X \leq x\},
\]
or the joint density function $\pi_{i}(x)=\frac{\partial \Pi_{i}(x)}{\partial x}$ for $x> 0$ and $i=0, 1, \ldots$, with $\pi_{i}(0)=\lim_{x\rightarrow 0^{+}}\pi_{i}(x)$.
According to the definition of the model, we have $\Pi_{i}(0)=0$ for $i\geq c$.
\bigskip

\textbf{Step~1:}
We can now define the transformations of interest. Let $\phi_{i}(\alpha)$ be the Laplace transform for $\pi_{i}(x)$, or
\[
\phi_{i}(\alpha)= \int_{0}^{\infty} \pi_{i}(x)e^{\alpha x}dx.
\]
Define $\psi(\alpha, z)=\sum_{i=c-1}^\infty \phi_{i}(\alpha) z^i$ and $\psi(z)=\sum_{i=c-1}^\infty \Pi_{i}(0)z^i$.
Then, the following fundamental form can be established.

\begin{theorem}[Theorem~1 in \cite{Li-Liu-Zhao:2019}] \label{the-phi-rel}
The fundamental equation can be rewritten as
 \begin{equation}\label{equ-funde}
   H(\alpha, z)\psi(\alpha, z)= \hat{H}_{1}(\alpha, z)\phi_{c-1}(\alpha)+  H_{2}(z)\psi(z)+ \hat{H}_{0}(\alpha, z),
 \end{equation}
 where
 \begin{align*}
   H(\alpha, z) &= -\lambda z^{2}+(-\alpha r+\lambda +c \mu) z- c\mu, \\
   \hat{H}_{1}(\alpha, z) &= \lambda z^{c}A_{c-2}(\alpha)+ H_{1}(\alpha, z), \\
   H_{2}(z)&=\lambda z^{2}-\lambda z-c \mu z+ c\mu,\\
    \hat{H}_{0}(\alpha, z) &= H_{0}(z)\Pi_{c-1}(0)+\lambda z^{c}\Pi_{c-2}(0)+\lambda z^{c}\sum_{n=0}^{c-2} \left[k_{n}\lambda^{c-2-n}\prod_{m=n}^{c-2}\frac{A_{m}(\alpha)}{(m+1)\mu}\right],
 \end{align*}
 with
 \[
     H_{1}(\alpha, z)=(\mu-\alpha r -\alpha)z^{c}-c\mu z^{c-1},
 \]
 \[
    H_{0}(z)=\mu z^{c}- c\mu z^{c-1},
 \]
 \[k_{0}=\mu\Pi_{1}(0)-\lambda\Pi_{0}(0),\]
\[
 k_{i}=\lambda\Pi_{i-1}(0)-(\lambda+i\mu)\Pi_{i}(0)+ (i+1)\mu \Pi_{i+1}(0),\ 1\leq i\leq c-2,
\]
and
\[
 A_{i}(\alpha)=\frac{(i+1)\mu}{\alpha+\lambda+i\mu-\lambda A_{i-1}(\alpha)},  \ 0\leq i\leq c-2, \ \ A_{-1}(\alpha)=0.
\]
\end{theorem}
We can see from the above that the kernel function $H(\alpha,z)$ is quadratic in $z$ (in the direction of the discrete state) and linear in $\alpha$ (in the direction of the continuous state), which is as expected from the previous discussions (quadratic for the discrete case and linear for the continuous case). We also point out that the boundary functions $ \hat{H}_{1}(\alpha, z)$ and $\hat{H}_{0}(\alpha, z)$ involved are more complicated than in the models discussed in the previous two sections, which place more challenges in the analysis, while the boundary function $H_{2}(z)$ is independent of $\alpha$.

It is helpful to know why we the specific fundamental equation in (\ref{equ-funde}) is adapted for the analysis. Recall that for the random walk in the quarter plane, there is only one layer of boundaries, the $x$-axis and the $y$-axis, and the fundamental form is the relationship between the joint distribution and two boundary distributions. For the  
fluid model driven by the $M/M/1$ queue, the situation is similar, but for the model driven by the $M/M/c$ queue, there are $c-1$ layers of boundaries. The fundamental equation in (\ref{equ-funde}) is a relationship between the joint distribution and the last layer $c-1$ of the boundary distribution. We also note that the presence of finite number of unknown in the fundamental form does not impact the analysis of the tail asymptotics. 
\bigskip

\textbf{Step~2:} Unlike in the previous sections, where the model is symmetric in the sense that the process of continuation for the two unknown functions is the same, the model discussed in this section is a combination of the discrete and continuous states. It requires separate treatments for the continuation, but the main required ingredients are the same: a discussion of the branch points, and the interplay between the two unknown functions.

For the kernel equation
\begin{equation}\label{equ-qua-form}
  H(\alpha, z) = az^{2}+ b(\alpha) z+ d =0,
\end{equation}
as a quadratic form in $z$,
where $a=-\lambda$, $b(\alpha)=-\alpha r+\lambda +c\mu$ and $d=- c\mu$, there are two solutions $z$ for a fixed $\alpha$, given by
\begin{equation}\label{equ-Z}
  Z_{\pm}(\alpha)=\frac{-b(\alpha)\pm \sqrt{\Delta(\alpha)}}{2a},
\end{equation}
$\Delta(\alpha)=b^{2}(\alpha)-4ad$ be the discriminant of the quadratic form; while for a fixed $z$, the only solution $\alpha$ is given by
\begin{equation}\label{equ-alpha-1}
 \alpha(z)=\frac{-\lambda z^{2}+ (\lambda+c\mu)z-c\mu}{zr}.
\end{equation}

This mixed model is special. Unlike for the pure discrete or pure continuous model (whose fundamental form is quadratic in both variables), the kernel function $H(\alpha, z)$ for the fluid model driven by the $M/M/c$ queue is quadratic in $z$ and linear in $\alpha$.  Therefore, for given $z$, $\alpha(z)$ is a single-valued function (see (\ref{equ-alpha-1})), while for given $\alpha$, a two-valued function $Z(\alpha)$ is determined. In this case, there are only two branch points (same as for the pure continuous case), whose explicit expressions are quite simple (see the next lemma).

\begin{lemma}[Lemma~1 in \cite{Li-Liu-Zhao:2019}] \label{lem-bran-point}
$\Delta(\alpha)$ has two positive zero points $\alpha_{1} = \frac{\left ( \sqrt{c\mu} - \sqrt{\lambda} \right )^2}{r}$ and $\alpha_{2} = \frac{\left ( \sqrt{c\mu} + \sqrt{\lambda} \right )^2}{r}$.
Moreover, $\Delta(\alpha)> 0$  in  $(-\infty, \alpha_{1})\cup (\alpha_{2}, \infty)$ and $\Delta(\alpha)< 0$  in $(\alpha_{1}, \alpha_{2})$.
\end{lemma}

The construction of the analytic branches of the algebraic function $Z(\alpha)$ is also simple. Similar to the previous discussions, the two branches $Z_0(\alpha)$ and $Z_1(\alpha)$ can be constructed by:
\[
  Z_{0}(\alpha)=Z_{-}(\alpha) \ \mbox{and} \ Z_{1}(\alpha)=Z_{+}(\alpha) \  \mbox{if} \  \Re(\alpha)> \frac{\lambda+c\mu}{r},
\]
\[
  Z_{0}(\alpha)=Z_{+}(\alpha) \ \mbox{and} \ Z_{1}(\alpha)=Z_{-}(\alpha) \  \mbox{if} \ \Re(\alpha) \leq \frac{\lambda+c\mu}{r}.
\]

Both $Z_0(\alpha)$ and $Z_1(\alpha)$ are analytic in the cut plane $\widetilde{\mathbb{C}}_{\alpha}= \mathbb{C}_\alpha \setminus \{[\alpha_{1},\alpha_{2}]\}$.
\begin{lemma}[Lemma~2 in \cite{Li-Liu-Zhao:2019}]\label{lem-ana-H}
The functions  $Z_{0}(\alpha)$ and $Z_{1}(\alpha)$ are analytic in $\widetilde{\mathbb{C}}_{\alpha}$.
Similarly, $\alpha(z)$
is meromorphic in $\mathbb{C}_z$ and $\alpha(z)$ has two zero points and one pole.
\end{lemma}
\bigskip

\textbf{Step~3:}
For the symmetric cases (states in both dimensions are discrete, or continuous, respectively) discussed in the previous sections, the analytic continuation process is also symmetric, in terms of the interplay of the two unknown functions. In the non-symmetric case considered in this section (states are discrete in one dimension and continuous in the other dimension), the interplay is not symmetric, for example $H_2(z)$ is a function of $z$ only, independent of $\alpha$, which makes the continuation of $\psi(z)$ much easier. However, since the function $\hat{H}_1(\alpha,z)$ (also $\hat{H}_0(\alpha,z)$) is much more complicated,  which makes the continuation for $\phi_{n-1}(\alpha)$ harder. The interplay relationship is given in the next lemma.

\begin{lemma}[Lemmas~6 and 7 in \cite{Li-Liu-Zhao:2019}] \label{lem-ana}
$\phi_{c-1}(\alpha)$ is analytic on $\{\alpha: \Re(\alpha) < \alpha_{dom}\}$, where $\alpha_{dom}=C_{p}(\phi_{c-1})> 0$ is the parameter of convergence for $\phi_{n-1}(\alpha)$;
and $\psi(z)$ is analytic on the disk $\Gamma_{z_{dom}}=\{z: |z|< z_{dom}\}$, where $z_{dom}=\frac{c\mu}{\lambda}$ is the radius of convergence for $\psi(z)$. Moreover,
$\phi_{c-1}(\alpha)$ can be analytically continued to the domain $D_{\alpha}=\{\alpha\in \widetilde{\mathbb{C}}_{\alpha}:  \hat{H}_{1}(\alpha, Z_{0}(\alpha))\neq 0\}\cap\{\alpha\in \widetilde{\mathbb{C}}_{\alpha}: |Z_{0}(\alpha) |< \frac{c\mu}{\lambda}\}$, and
\begin{equation}\label{equ-asy-ana-1}
   \phi_{c-1}(\alpha)=-\frac{H_{2}( Z_{0}(\alpha))\psi(Z_{0}(\alpha))+ \hat{H}_{0}(\alpha, Z_{0}(\alpha))}{ \hat{H}_{1}(\alpha, Z_{0}(\alpha))};
\end{equation}
and
$\psi(z)$ can be analytically continued to the domain $D_{z}=\{z\in \mathbb{C}: H_{2}(z)\neq 0\}\cap\{z\in \mathbb{C}: \Re(\alpha(z)) < \alpha_{dom}\}$ and
\begin{equation}\label{equ-asy-ana-2}
   \psi(z) =-\frac{\hat{H}_{1}(\alpha(z), z) \phi_{c-1}(\alpha(z))+\hat{H}_{0}(\alpha(z), z)}{ H_{2}(z)}.
\end{equation}
\end{lemma}

In order to specify detailed asymptotic property at the dominant singularity, we need to classify the type of the dominant singularity. For this purpose and based on the above lemma, we can have:
\begin{lemma}[Lemma~8 in \cite{Li-Liu-Zhao:2019}] \label{lem-pol-ana}
The convergence parameter $\alpha_{dom}$ satisfies $0< \alpha_{dom} \leq \alpha_{1}$. If $\alpha_{dom} < \alpha_{1}$, then $\alpha_{dom}=\alpha^*$ is necessarily a zero point of $\hat{H}_{1}(\alpha, Z_{0}(\alpha))$.
\end{lemma}
It is clear now that there are only two candidates for the dominant singularity for $\phi_{n-1}(\alpha)$, which are $\alpha_1$, the leftmost branch point greater than 0, and $\alpha^*$, the leftmost zero for $\hat{H}_{1}(\alpha, Z_{0}(\alpha))$ in $(0,\alpha_1]$. Therefore, there are a total of three cases for the asymptotic property.
For convenience, let $\alpha^*>\alpha_1$ if there is no zero of $\hat{H}_{1}(\alpha, Z_{0}(\alpha))$ in $(0,\alpha_1]$. Then, the three cases about the dominant singularity of $\phi_{n-1}(\alpha)$ are given by:
\begin{description}
\item[Case~1.] $\alpha_{dom} = \alpha^* < \alpha_{1}$;

\item[Case~2.] $\alpha_{dom} = \alpha^* = \alpha_{1}$;

\item[Case~3.] $\alpha_{dom} = \alpha_{1} < \alpha^*$.
\end{description}

For convenience, we make the following assumptions:
\begin{assumption}\label{ass-1}
\begin{description}
\item[(i)] The function $\hat{H}_{1}(\alpha, Z_{0}(\alpha))$ has at most one real zero point in $(0, \alpha_{1}]$, denoted  by $\tilde{\alpha}$, if such a zero exists.

\item[(ii)] The zero point $\tilde{\alpha}$ satisfies  $H_{2}( Z_{0}(\tilde{\alpha}))\psi(Z_{0}(\tilde{\alpha}))+ \hat{H}_{0}(\tilde{\alpha}, Z_{0}(\tilde{\alpha}))\neq 0$.

\item[(iii)] The unique $\tilde{\alpha}$ is a multiple zero of $k$ times for function $\hat{H}_{1}(\alpha, Z_{0}(\alpha))$, where $k \geq 1$ is an integer.
\end{description}
\end{assumption}
These conditions can be justified for some models, and we conjectured that these conditions are automatically satisfied by the fluid model driven by the $M/M/c$ queue for a general $c$. See more discussion in \cite{Li-Liu-Zhao:2019}. The above assumptions are necessary when asymptotic analysis is performed at the dominant singularity. For example,  without these assumptions, we cannot guarantee the constant $g$ in the Tauberian-like theorem (Theorem~\ref{tauberian-1}) is non-zero. This is a common challenge for functions with a complex structure (see also discussion in Section~\ref{sec:6}.

The detailed asymptotic property at the dominant singularity is summarized in the following theorem.
\begin{theorem}[Theorems~2 and 3 in \cite{Li-Liu-Zhao:2019}] \label{the-tail-asy-Pi1}
\textbf{1.}
Suppose that Assumption~\ref{ass-1} is satisfied.
Then, for the function $\phi_{c-1}(\alpha)$, a total of three types of asymptotics exist, corresponding to the above three cases, as $\alpha$ approaches to $\alpha_{dom}$:
\begin{description}
\item[Case~1.]
\[
   \lim_{\alpha\rightarrow \alpha_{dom}=\alpha^*}(\alpha_{dom} -\alpha)^{k}\phi_{c-1}(\alpha)= c_{1},
\]
where $k$ is given is Assumption~\ref{ass-1};

\item[Case~2.]
\[
   \lim_{\alpha\rightarrow \alpha_{dom}=\alpha^*=\alpha_1}\sqrt{\alpha_{dom}-\alpha}\cdot\phi_{c-1}(\alpha)=c_{2};
\]

\item[Case~3.]
\[
   \lim_{\alpha\rightarrow \alpha_{dom}=\alpha_1}\sqrt{\alpha_{dom}-\alpha}\cdot\phi_{c-1}'(\alpha)=c_{3},
\]
\end{description}
where
$c_i$, $i=1,2,3$, are non-zero constants, whose expressions can be found in Theorem~2 of \cite{Li-Liu-Zhao:2019}.
\medskip

\textbf{2.} For the function $\psi(z)$, we have the following asymptotic property as
$z$ approaches to $z_{dom}=\frac{c\mu}{\lambda}$:
\[
    \lim_{z\rightarrow z_{dom}}(z_{dom}-z)\psi(z)= d,
\]
where $d$ is a non-zero constant, whose expression can be found in Theorem~3 of \cite{Li-Liu-Zhao:2019}.
\end{theorem}
\bigskip

\textbf{Step~4:}
In the final step, we apply an Tauberian-like theorem, either Theorem~\ref{5-thm1} (along the $\alpha$ direction) or Theorem~\ref{tauberian-1} (along the $z$ direction) to derive the exact tail asymptotic properties of $\pi_{c-1}(x)$ and $\Pi_i(0)$, respectively.

\begin{lemma}[Lemmas~9 and 10 in \cite{Li-Liu-Zhao:2019}] \label{lem-tail-asy-2}
\textbf{1.}  Under Assumption~\ref{ass-1}, the density function $\pi_{c-1}(x)$ of the fluid queue has the following tail asymptotic properties corresponding to the three cases:
\begin{description}
\item[Case~1.]
\[
   \pi_{c-1}(x)\sim C_{1}e^{-\alpha_{dom}\, x}x^{k-1},
\]
where $k$ is given is Assumption~\ref{ass-1};

\item[Case~2.]
\[
   \pi_{c-1}(x) \sim C_{2} e^{-\alpha_{dom}\, x}x^{-\frac{1}{2}};
\]

\item[Case~3.]
\[
  \pi_{c-1}(x) \sim C_{3} e^{-\alpha_{dom}\, x}x^{-\frac{3}{2}},
\]
\end{description}
where $C_{1}=\frac{c_{1}}{\Gamma(k)}$, $C_{2}=\frac{c_{2}}{\sqrt{\pi}}$ and $C_{3}=\frac{-c_{3}}{2\sqrt{\pi}}$ with $c_{i}, i=1, 2, 3$, being given in Theorem~\ref{the-tail-asy-Pi1}.

\textbf{2.}
The boundary probabilities $\Pi_{i}(0)$ of the fluid queue has the following tail asymptotic property:
 \[
      \Pi_{i}(0)\sim d \Big (\frac{1}{z_{dom}} \Big )^{i+1},
 \]
where $z_{dom}=\frac{c\mu}{\lambda}$ and $d$ is given in Theorem~\ref{the-tail-asy-Pi1}.
\end{lemma}

The tail asymptotic property for $\pi_i(x)$ and $\Pi_i(x)$ (for $i>c-1$) can be obtained through recursive formulas, which is stated in the next theorem. 
The recursive scheme used for this model is also valid for other models, for example the models discussed in Sections~\ref{sec:3} and \ref{sec:4}. 
\begin{theorem}[Theorem~5 in \cite{Li-Liu-Zhao:2019}]  \label{the-tail-asy-4}
Under Assumption~\ref{ass-1}, the buffer content (or the marginal) distribution $\Pi(x)$ of the fluid queue has the following tail asymptotic properties according to the three cases:

\begin{description}
\item[Case~1.]
\[
   \Pi(x)- 1 \sim -\frac{\tilde{C}_{1}}{\alpha_{dom}} e^{-\alpha_{dom}\, x}x^{k-1},
\]
where $k$ is given is Assumption~\ref{ass-1};

\item[Case~2.]
\[
    \Pi(x)- 1 \sim -\frac{\tilde{C}_{2}}{\alpha_{dom}} e^{-\alpha_{dom}\, x}x^{-\frac{1}{2}};
\]

\item[Case~3.]
\[
   \Pi(x)-  1 \sim -\frac{\tilde{C}_{3}}{\alpha_{dom}}e^{-\alpha_{dom}\, x}x^{-\frac{3}{2}},
\]
\end{description}
where
  expressions for $\tilde{C}_{i}$, $i=1, 2, 3$, can be found in Theorem~5 in \cite{Li-Liu-Zhao:2019}.
\end{theorem}

\begin{example}[Driven by the $M/M/1$ queue]
For this example, we can show that $\hat{H}_{1}(\alpha, Z_{0}(\alpha))=0$ has a unique zero $\alpha^*=\frac{\mu}{r+1}-\lambda \leq \alpha_1$ in $(0,\alpha_1]$ such that
$H_{2}( Z_{0}(\alpha^*))\psi(Z_{0}(\alpha^*))+ \hat{H}_{0}(\alpha^*, Z_{0}(\alpha^*))\neq 0$.
 This implies that Assumption~\ref{ass-1} holds and only two asymptotic cases exist:

\textbf{Case~1.} If $\alpha^*=\frac{\mu}{r+1}-\lambda < \alpha_{1} = \frac{\left ( \sqrt{\mu} - \sqrt{\lambda} \right )^2}{r}$, then
\[
  \Pi(x)- 1 \sim -\frac{(r+1) \alpha^*}{r}  c_{1} e^{-\tilde{\alpha} x};
\]

\textbf{Case~2.} If $\alpha^*=\frac{\mu}{r+1}-\lambda = \alpha_{1} = \frac{\left ( \sqrt{\mu} - \sqrt{\lambda} \right )^2}{r}$, then
\[
      \Pi(x)- 1 \sim -\frac{(r+1)\alpha^*}{r\sqrt{\pi}} c_{2} e^{-\alpha^* x}x^{-\frac{1}{2}}.
\]
Here  $c_{1}$ and $c_{2}$ are defined in Theorem~\ref{the-tail-asy-Pi1}.

\end{example}

\section{Possible extensions and discussions} \label{sec:6}

In this section, we briefly discuss the possibility of extending the kernel method to generalized problems, mainly focusing on the new challenges.

\subsection{Absorbing time to boundaries of the random walk in the quarter plane}

In this subsection, we demonstrate that the kernel method discussed in this paper can be extended to deal with the absorbing time to the boundaries for a random walk in the quarter plane. To be specific, we consider a voter model, studied by Kurkova and Raschel~\cite{Kurkova-Raschel:2013}, in which exact tail asymptotic properties for the absorbing time were obtained based on the solution of the unknown functions for the absorbing time. In terms of the kernel method, there is no need to have a complete characterization of the unknown function. The same asymptotic results were obtained (see, Li, Song and Zhao~\cite{Li-Song-Zhao:2014}.

The voter model is a two-dimensional Markov chain, or a random walk in the quarter plane, $(X(k), Y(k))$ absorbed to a boundary with probability one.  For a detailed model description, we refer readers to \cite{Kurkova-Raschel:2013}. We are interested in the absorbing time:
\[
    \tau =\inf \{k\in \mathbb{Z}_{+}:X(k)=0\;\text{ or }\;Y(k)=0\},
\]
which is equal to the minimal absorbing times, $\tau =\min \{S,T\}$, to the horizontal and vertical axes, respectively, defined by:
\begin{equation*}
S=\inf \{k\in \mathbb{Z}_{+}:Y(k)=0\},\;\;\;T=\inf \{k\in \mathbb{Z}
_{+}:X(k)=0\}.
\end{equation*}
By symmetry, we only provide discussions for the tail asymptotic of the conditional probabilities $P_{(i_{0},j_{0})}[S=k]$, absorbed to the horizontal axis at time $k$ given that the initial state $(i_{0},j_{0})$ is not on either axis (boundary).

Let
\[
    \pi(i,j)=\pi_{(i_{0},j_{0})}(i,j)= P[(X(k),Y(k))=(i,j)|(X(0),Y(0))=(i_{0},j_{0})].
\]
The fundamental form for this problem can be form as:
\begin{equation*}
    K(x,y,z)H^{i_{0},j_{0}}(x,y,z)=h^{i_{0},j_{0}}(x,z)+\tilde{h}^{i_{0},j_{0}}(y,z)-x^{i_{0}}y^{j_{0}},
\end{equation*}
where
\begin{align*}
    K(x,y;z) &= xyz\left[ \sum_{-1\leq i,j\leq 1}p_{i,j}x^{i}y^{j}-1/z\right], \\
    H^{i_{0},j_{0}}(x,y;z) &= \sum_{i>1}\sum_{j>1}\sum_{k \geq 0} \pi(i,j) x^{i-1}y^{j-1}z^{k}, \\
    h^{i_{0},j_{0}}(x;z) &= \sum_{i>1}\sum_{k\geq 0} \pi(i,0) x^{i}z^{k}, \\
    \tilde{h}^{i_{0},j_{0}}(y;z) &= \sum_{j>1}\sum_{k\geq 0} \pi(0,j) y^{j}z^{k}.
\end{align*}
The new challenge here is the present of the third (time) variable in the fundamental form. However, for our purpose, we only need to consider the functions
$h^{i_0, j_0}(1;z)$ and $\tilde{h}^{i_0, j_0}(y; z)$. In this way, the problem is converted to a two-dimensional one, for which the kernel method can prevail.
The key using the kernel method is to show the asymptotic property of
the function $h^{i_{0},j_{0}}(1;z)$ at $z=1$:
\begin{equation*}
\lim_{z\rightarrow 1}(1-z)^{-3/2}h^{i_{0},j_{0}}(1;z)=g,
\end{equation*}%
or
\begin{equation*}
\lim_{z\rightarrow 1}(1-z)^{-1/2}\frac{dh^{i_{0},j_{0}}(1;z)}{dz}=c.
\end{equation*}

Notice that the absorbing probability $P_{(i_{0},j_{0})}[S=k]$ into $x$-axis can expressed by
\[
    P_{(i_{0},j_{0})}[S=k] = \sum_{i} \pi(i,0),
\]
which is the coefficient of $z^k$ of the function $h^{i_{0},j_{0}}(1;z)$. Applying the Tauberian-like theorem leads to:
\[
    P_{(i_{0},j_{0})}[S=k] \sim c k^{-5/2}.
\]

\subsection{Random walks in the quarter plane modulated by a finite-state Markov chain}

The random walk in the quarter plane modulated by a Markov chain is a generalization fo the random walk in the quarter plane with the state space $S=\{(m,n; l): m, n = 0, 1, \ldots; l = 1, 2, \ldots, M \}$ and transition diagram given in Figure~1, where $A_{i,j}$ and $A^{(k)}_{i,j}$ are all nonnegative matrices of size $M \times M$ with $M < \infty$, and $\sum_{i,j} A_{i,j}$ and $\sum_{i,j} A^{(k)}_{i,j}$ ($k =0, 1, 2$) are stochastic.

\begin{figure}[h] \label{fig:tandem}
\centering
\includegraphics[width=3.8in]{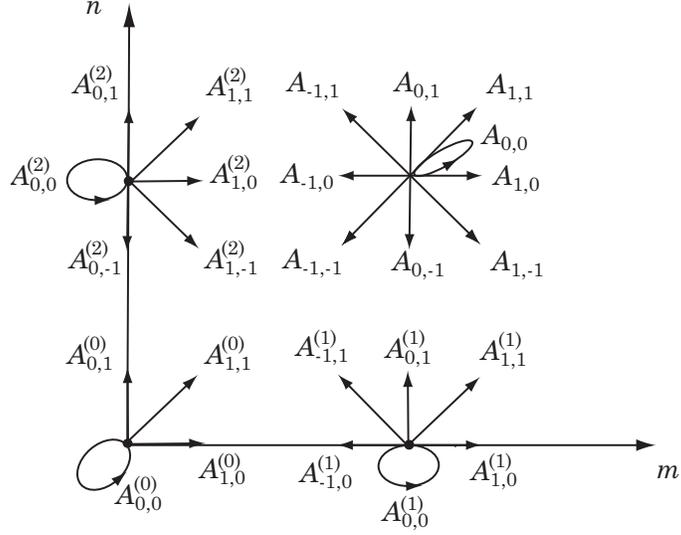}
\caption{Transition diagram for Markov modulated random walks}
\end{figure}

Under a stability condition (see, for example \cite{Ozawa-Kobayashi:2018}), let $\pi_{m,n;k}$  be the stationary probability vector. Then, the fundamental form (in matrix form) is given by:
\begin{equation} \label{eqn:FF}
    -\Pi(x,y) H(x,y) = \Pi_1(x) H_1(x,y)  + \Pi_2(y) H_2(x,y)  + \Pi_0 H_0(x,y),
\end{equation}
where
\begin{align*}
         H(x,y) =& xy \left (I-\sum^1_{i=-1}\sum^1_{j=-1}x^iy^jA_{ij} \right ), \\
         H_1(x,y) =& x \left (I-\sum^1_{i=-1}\sum^1_{j=0}x^iy^jA_{ij}^{(1)} \right ), \\
         H_2(x,y) =& y \left (I-\sum^1_{i=0}\sum^1_{j=-1}x^iy^jA_{ij}^{(2)} \right ),\\
         H_0(x,y) =& \left (I-\sum^1_{i=0}\sum^1_{j=0}x^iy^jA_{ij}^{(0)} \right ), 
\end{align*}
\begin{align*}
         \Pi(x,y) =& \left (\sum^\infty_{i=1}\sum^\infty_{j=1}\pi_{i,j;1}x^{i-1}y^{j-1},\sum^\infty_{i=1}\sum^\infty_{j=1}\pi_{i,j;2}x^{i-1}y^{j-1},\ldots, \sum^\infty_{i=1}\sum^\infty_{j=1}\pi_{i,j;M}x^{i-1}y^{j-1} \right )_{1\times M}, \\
         \Pi_1(x) =&  \left (\sum^\infty_{i=1}\pi_{i,0;1}x^{i-1},\sum^\infty_{i=1}\pi_{i,0;2}x^{i-1},\ldots,\sum^\infty_{i=1}\pi_{i,0;M}x^{i-1}\right )_{1\times M},\\
         \Pi_2(y) =&  \left (\sum^\infty_{j=1}\pi_{0,j;1}y^{j-1},\sum^\infty_{j=1}\pi_{0,j;2}y^{j-1},\ldots,\sum^\infty_{j=1}\pi_{0,j;M}y^{j-1}\right )_{1\times M},\\
         \Pi_0 =& \left (\pi_{0,0;1},\pi_{0,0;2},\ldots,\pi_{0,0;M}\right )_{1\times M}.\\
\end{align*}

A special example was considered in Liu, Wang and Zhao~\cite{Liu-Wang-Zhao:2016}, which is a Markov modulated two-demand queueing system. Specifically, 
\[
    A_{1,1}= \left [ \begin{array}{cc} \lambda_0 p & \lambda_0 \bar{p} \\ \lambda_1 \bar{q} & \lambda_1 q \end{array} \right ], \quad
A_{0,-1}= \mu_2 \left [ \begin{array}{cc} p & \bar{p} \\ \bar{q} & q \end{array} \right ],
\]
\[
A_{-1,0}= \mu_1 \left [ \begin{array}{cc} p & \bar{p} \\ \bar{q} & q \end{array} \right ], \quad
A_{0,0}= \left [ \begin{array}{cc} \lambda_1 & 0 \\ 0 & \lambda_0 \end{array} \right ],
\]
and
\[
    A_{0,0}^{(0)}=\left [ \begin{array}{cc} 1-\lambda_0 & 0 \\ 0 & 1-\lambda_1 \end{array} \right ], \quad
    A_{0,0}^{(1)}=\left [ \begin{array}{cc} \lambda_1+\mu_2 & 0 \\ 0 & \lambda_0+\mu_2 \end{array} \right ],
\]
\[
    A_{0,0}^{(2)}=\left [ \begin{array}{cc} \lambda_1+\mu_1 & 0 \\ 0 & \lambda_0+\mu_1 \end{array} \right ],
\]
$A_{1,1}^{(0)}=A_{1,1}^{(1)}=A_{1,1}^{(2)}=A_{1,1}$, $A_{-1,0}^{(1)}=A_{-1,0}$ and $A_{0,-1}^{(2)}=A_{0,-1}$.
The key idea to extend the kernel method for this type of model is to convert the fundamental equation of matrix-form into one of scaler-form, such that the standard kernel method can be applied. In \cite{Liu-Wang-Zhao:2016}, the maximal eigenvalue of a matrix was used for the conversion. However, their idea could is valid only for some special cases. A method combining the maximal eigenvalue with the conjugate of a matrix could be a promising one for the general model, which will be addressed in a separate work.

The random walk in the quarter plane modulated by a Markov chain was also considered in \cite{Ozawa-Kobayashi:2018}, which was referred to as the discrete-time, two-dimensional quasi-birth-and-death process (2d-QBD). Instead of special examples, the authors of \cite{Ozawa-Kobayashi:2018} considered a general setting for this type of Markovian model. Their method can be considered as a combination of the geometric method, introduced by Miyazawa in \cite{Miyazawa:08}, and the asymptotic analysis, an important component of the kernel method. The exact tail behaviour was characterized based on the same key ideas as that of the kernel method: the property of the kernel function in block form, the dimension reduction based on the interplay between two unknown functions in vector form, and the exact tail asymptotic property as the consequence of applying the Tauberian-like theorem.
The additional challenge for this type of models is due to the presence of the third finite dimension, which leads to the matrix-form fundamental equation (\ref{eqn:FF}) and two vector-form unknown generating functions $\Pi_1$ and $\Pi_2$. To meet this additional challenge, spectral analysis in the matrix-analytic method is the bridge to connect matrices to scalers, and a $G$ matrix ($G_1(z)$ in \cite{Ozawa-Kobayashi:2018}) was introduced for the interplay of the two unknown functions. However, due to the technical difficulties, the constant $g$ given in the Tauberian-like theorem could not be proved to be non-negative in general, which is expected to be true. Therefore, the constants $l_1$ and $l_2$ in the main theorem (Theorem~2.1 in \cite{Ozawa-Kobayashi:2018}) could not be confirmed to be 1 that would ensure that the 2d-QBD has the same types of tail asymptotics as that for the random walk in the quarter plane.

\subsection{Higher-dimensional problems}

There are two main challenges for extending the kernel method for higher-dimensional random walks: (1) analytic continuation of the unknown functions; and (2) additional information required about relationships between the unknown functions. Although the key idea for analytic continuation of the unknown functions for two-dimensional problems is often still valid for higher-dimensional problems, the continuation process requires much more technical details. The right-hand side of a fundamental form for a $k$-dimensional model involves more (up to $2^k -2$) unknown functions of up to $k-1$ variables. More information (or equations) about these unknown functions need to be revealed based on the fundamental form. Until now, we have not seen any systematic treatment for $k$-dimensional models. However, for specific cases, one may find ways to reduce the dimension such that the kernel method for two-dimensional problems can be applied.

One of such moelds was considered in Dai, Dawson and Zhao~\cite{Dai-Dawson-Zhao:2018}, which is a three-dimensional Brownian-driven tandem queue with intermediate inputs.
The buffer content $L_i(t)$ at node $i$ at time $t\geq 0$, for $i=1,2,3$, is described by
\begin{align*}
    L_1(t)&=L_1(0)+X_1(t)-c_1t+Y_1(t), \\
    L_i(t)&=L_i(0)+X_i(t)+c_{i-1}t-c_{i} t-Y_{i-1}(t)+Y_i(t),\; i\geq 2,
\end{align*}
where $X_i(t)$ is the exogenous input to node $i$, assumed to be a Brownian process of the form:
\[
 X_i(t)=\lambda_i t+B_i(t),
\]
with $\lambda_i>0$ and $B_i(t)$ being a Brownian motion with variance $\sigma^2_i$ and no drift,
and $Y_i(t)$ is a regulator at node $i$, which is a minimal nondecreasing process for $L_i(t)$ to be nonnegative.
Under a stability condition (for example, see \cite{Dai-Dawson-Zhao:2018}), let $(L_1,L_2,L_3)'$ be the stationary vector of $(L_1(t),L_2(t),L_3(t))$ having the same distribution $\pi$ as the stationary distribution of the model, we are interested in the exact tail asymptotic of $P(L_3 >t)$.
Notice that the tail asymptotic property for $P(L_1 >t)$, or $P(L_2 >t)$ is trivial, or can be obtained through the kernel method directly, respectively.
This is a continuous-state model, we can similarly (to what has been done in Section~\ref{sec:4}) define the boundary measures
$V_i(\cdot)$,  $i=1,2,3$, by
\[
    V_i(A)=\E_\pi\Big[\int_0^1 I_{\{L(u)\in A\}}d Y_i(u)\Big].
\]
We can then establish the fundamental form for the model:
\[
  -H(x,y,z)\phi(x,y,z)=H_1(x,y)\phi_1(x,y,z)+H_2(y,z)\phi_2(x,y,z)+H_3(z)\phi_3(x,y,z),
\]
where
\begin{align*}
    \phi(x,y,z)&=\E \big[e^{ x L_1+ y L_2+ z L_3}\big], \\
    \phi_i(x,y,z)&=\int_{\R_+^3}e^{<w,\theta> }V_i(d\theta)=\E\Big[\int_0^1 e^{<w,\;L(t)>}dY_i\Big], \quad i=1,2,3,
\end{align*}
are unknown functions, and
\begin{align*}
    H(x,y,z)&= \frac{1}{2}\big(\sigma_1^2x^2+\sigma_2^2y^2 +\sigma_3^2z^2\big) -(c_1-\lambda_1)x -(c_2-\lambda_2-c_{1})y -(c_3-\lambda_3-c_{2})z,\\
    H_1(x,y)&=x-y, \\
    H_2(y,z)&=y-z, \\
    H_3(z)&=z.
\end{align*}

For this specific model, the analytic continuation of the unknown functions $\phi_i$, $i=1,2, 3$, and their asymptotic property at the dominant singularity were obtained in \cite{Dai-Dawson-Zhao:2018} based on relationships (or interplay) of the unknown functions. The key idea for the continuation is equivalent to a dimension reduction procedure. This procedure is quite complicated, in which many technical details were involved, and unfortunately it could not be extended to a proof for the case with dimensions bigger than three. The asymptotic analysis is also not simple, but the main idea is also valid for higher-dimensional tandem queues. The application of the Tauberian-like theorem remains straightforward to give the following tail asymptotic result:
\begin{theorem}[Theorem~4.1 in \cite{Dai-Dawson-Zhao:2018}] \label{C-1}
For the marginal distribution  $\P(L_3>t)$, we have the following tail asymptotic properties:
\begin{description}
\item[Case~1.] If $z_{dom}=z^*<z^{max}$,
then
\[
    \P(L_3>t) \sim \hat{K}_1e^{-z_{dom} t};
\]

\item[Case~2.] If $z_{dom}=z^{max}<z^*$, then
\[
    \P(L_3>t) \sim \hat{K}_2 e^{-z_{dom} t} t^{-\frac{3}{2}};
\]

\item[Case~3.] If $z_{dom}=z^*=z^{max}$, then
\[
    \P(L_3>t) \sim \hat{K}_3 e^{-z_{dom} t} t^{-\frac{1}{2}},
\]
\end{description}
where both $z^*$ and $z^{max}$ are explicitly given in \cite{Dai-Dawson-Zhao:2018}, and
$\hat{K}_i$, $i=1,2,3$, are non-zero constants.
\end{theorem}

\vspace*{3mm}
\noindent \textbf{Acknowledgements:} The author thank the anonymous reviewer and the guest editor for their constructive valuable comments and suggestions for the improvement of the presentation of this work, and acknowledges that this work was supported in part through a Discovery Grant of NSERC.


\begin{thebibliography}{99}

\bibitem{Adan-Foley-McDonald:09} Adan, I., Foley, R.D. and McDonald,
D.R. (2009) Exact asymptotics for the stationary distribution of a
Markov chain: a production model, \textit{Queueing Systems},
\textbf{62}, 311--344.

\bibitem{AWZ:1993} Adan, I.J.B.F., Wessels, J. and Zijm, W.H.M. (1993)
A compensation approach for two-dimensional Markov processes,
\textit{Advances in Applied Probability},  \textbf{25}, 783--817.

\bibitem{ANY:2014}
Avrachenkov, K., Nain, P. and Yechiali, U. (2014)
A retrial system with two input streams and two orbit queues, 
\textit{Queueing Systems}, \textbf{77}(1), 1--`31.

\bibitem{B-BM-D-F-G-GB:02}
Banderier, C., Bousquet-M\'{e}lou, M,  Denise, A, Flajolet, P.,
Gardy, D. and Gouyou-Beauchamps, D. (2002) Generating functions of
generating trees, \textit{Discrete Math.}, \textbf{246}, 29--55.

\bibitem{Bender:1974} Bender, E. (1974)
Asymptotic methods in enumeration, \textit{SIAM Review},
\textbf{16}, 485--513.

\bibitem{Borovkov-M:01}
Borovkov, A.A. and Mogul'skii, A.A. (2001) Large deviations for
Markov chains in the positive quadrant, \textit{Russian Math.
Surveys}, \textbf{56}, 803--916.

\bibitem{Bousquet-Melou:05} Bousquet-M\'{e}lou, M. (2005) Walks in the quarter plane: Kreweras?algebraic model, \textit{Annals
of Applied Probability}, \textbf{15}, 1451-?491.

\bibitem{cohen-boxma:83} Cohen, J. W. and Boxma, O. J. (1983) \textit{Boundary Value
Problems in Queueing System Analysis}, North-Holland, Amsterdam.

\bibitem{Dai-Dawson-Zhao:2015} Dai, H., Dawson, D.A. and Zhao, Y.Q. (2015)
Kernel method for staionary tails: from discrete to continuous, in \textit{Asymptotic Laws and Methods in Stochastics}, edited by D. Dawson, R. Kulik, M. Ould Haye, B. Szyszkowicz and Y. Zhao, 297--328, Springer, New York.

\bibitem{Dai-Dawson-Zhao:2018}
Dai, H., Dawson, D. and Zhao, Y.Q. (2018)
Exact tail asymptotics for a three dimensional Brownian
tandem queue with intermediate inputs, under revisions. (https://arxiv.org/abs/1807.08425)

\bibitem{Dai-Kong-Song:2017} Dai, H., Kong, L. and Song, Y. (2017)
Exact tail asymptotics for a two-stage queue: Complete solution via kernel method, \textit{RAIRO-Oper. Res.}, \textbf{51}, 1211--1250.

\bibitem{Dai-Zhao:2013} Dai, H. and Zhao, Y.Q. (2013) Wireless 3-hop networks with stealing revisited: A kernel approach. \textit{INFOR}, \textbf{51(4)}, 192-205.

\bibitem{DH1992}
Dai, J.G. and Harrison, J.M. (1992)
Reflected Brownian motion in an orthant: numerical methods for steady-state analysis,
\textit{Ann. Appl. Probab.}, \textbf{2}, 65--86.

\bibitem{DM2011}
Dai, J.G. and Miyazawa, M. (2011)
Reflecting Brownian motion in two dimensions: Exact asymptotics for the stationary distribution,
\textit{Stochastic systems}, \textbf{1}, 146--208.

\bibitem{DM2012}
Dai, J.G. and Miyazawa, M. (2013)
Stationary distribution of a two-dimensional SRBM: geometric views and boundary measures,
\textit{Queueng Syst.}, \textbf{74}, 181--217.

\bibitem{FI:1979} Fayolle, G. and Iasnogorodski, R.  (1979)
Two coupled processors: the reduction to a
Riemann-Hilbert problem, \textit{Z. Wahrscheinlichkeitsth}, \textbf{47}, 325--351.

\bibitem{FIM:2017} Fayolle, G., Iasnogorodski, R. and Malyshev, V. (2017)
\textit{Random Walks in the Quarter-Plane}, second ed., Springer, New York.

\bibitem{FKM:82}
Fayolle, G., King, P.J.B. and  Mitrani, I. (1982)
The solution of certain two-dimensional Markov models,
 \textit{Adv. Appl. Prob.}, \textbf{14}, 295--308.


\bibitem{Feller:1971} Feller, W. (1971) \textit{An Introductioin to Probability Theory and Its Applications}, Vol. II, second ed., Wiley,
New York.

\bibitem{FO:1990} Flajolet, P. and Odlyzko, A. (1990)
Singularity analysis of generating functions, \textit{SIAM J.
Disc. Math.}, \textbf{3}, 216--240.


\bibitem{Flatto-McKean:77}
Flatto, L. and McKean, H.P. (1977)
Two queues in parallel, \textit{Comm. Pure Appl. Math.}, \textbf{30}, 255--263.

\bibitem{Flatto-Hahn:84}
Flatto, L. and Hahn, S. (1984)
Two parallel queues created by arrivals with two demands I,
\textit{SIAM J. Appl. Math.}, \textbf{44}, 1041--1053.

\bibitem{Flatto:85}
Flatto, L. (1985)
Two parallel queues created by arrivals with two demands II,
\textit{SIAM J. Appl. Math.}, \textbf{45}, 861--878.


\bibitem{Flajolet-Sedgewick:09} Flajolet, F. and Sedgewick, R.
(2009) \textit{Analytic Combinatorics}, Cambridge University
Press.

\bibitem{Foley-McDonald:01} Foley, R.D. and McDonald, D.R. (2001)
Join the shortest queue: stability and exact asymptotics, \textit{Annals of Applied Probability}, \textbf{11}, 569--607.

\bibitem{Foley-McDonald:05a} Foley, R.D. and McDonald, R.D. (2005) Large deviations of a modified Jackson network:
stability and rough asymptotics, \textit{Annals of Applied Probability}, \textbf{15}, 519-541.

\bibitem{Foley-McDonald:05b} Foley, R.D. and McDonald, R.D. (2005)
Bridges and networks: exact asymptotics, \textit{Annals of Applied Probability}, \textbf{15}, 542--586.

\bibitem{GKL2011}
Guillemin, F., Knessl, C. and Van Leeuwaarden, J.S.H. (2011)
Wireless multi-hop networks with stealing: large buffer asymptotics via the Ray method,
\textit{SIAM Journal on Applied Mathemtics},  \textbf{71}, 1220--1240.

\bibitem{Guillemin-Leeuwaarden:09} Guillemin, F. and van Leeuwarden, J.S.H. (2011)
Rare event asymptotics for a random walk in the quarter plane,
\textit{Queueing Systems}, \textbf{67(1)}, 1--32.

\bibitem{Guillemin-Pinchon:2004}
Guillemin, F. and Pinchon, D. (2004)
Analysis of generalized processor-sharing systems with two classes of customers and exponential services,
\textit{Journal of Applied Probability}, \textbf{41(03)}, 832--858.

\bibitem{Haque:03} Haque, L. (2003)
textit{Tail Behaviour for Stationary Distributions for Two-Dimensional Stochastic Models}, Ph.D. Thesis, Carleton
University, Ottawa, ON, Canada.

\bibitem{Haque-Liu-Zhao:05}
Haque, L., Liu, L. and Zhao, Y.Q. (2005) Sufficient conditions for a
geometric tail in a QBD process with countably many levels and
phases, \textit{Stochastic Models}, \textbf{21}(1), 77--99.

\bibitem{Harrison-Williams:1987}
Harrison, J.M. and Williams, R.J. (1987) Brownian models of open queueing networks with homogeneous customer populations, \textit{Stochastics}, \textbf{22}(2), 77--115.

\bibitem{He-Li-Zhao:08} He, Q., Li, H. and Zhao, Y.Q. (2009)
Light-tailed behaviour in QBD process with countably many phases,
\textit{Stochastic Models}, \textbf{25}, 50--75.

\bibitem{Houtum-etc:2001} Houtum, G.J. van, Adan, I.J.B.F., Wessels,
J. and Zijm, W.H.M. (2001) Performance analysis of parallel
identical machines with a generalized shortest queue arrival
mechanism, {\it OR Spektrum}, \textbf{23}, 411--428.

\bibitem{Khanchi:08} Khanchi, Aziz (2008) \textit{State of a network when one node overloads},
Ph.D. Thesis, University of Ottawa.

\bibitem{Khanchi:09} Khanchi, Aziz (2009) Asymptotic hitting distribution for a reflected random walk in the
positive quardrant, \textit{Stochastic Models}, \textbf{27}, 169--201.

\bibitem{Kingman:1961} Kingman, J.F.C. (1961) Two similar queues in parallel, \textit{Ann. Math. Statist}. \textbf{32}, 1314--1323.

\bibitem{Knuth:69} Knuth, D.E. (1969)
\textit{The Art of Computer Programming, Fundamental Algorithms},
vol. 1, second ed., Addison-Wesley.

\bibitem{Kobayashi-Miyazawa:2013} Kobayashi, M. and  Miyazawa, M. (2013)
Revisiting the tail asymptotics of the double QBD process: Refinement and complete solutions for the coordinate and diagonal directions,
in \textit{Matrix-Analytic Methods in Stochastic Models}, edited by G. Latouche, V. Ramaswami, J. Sethuraman,  K. Sigman,
M.S. Squillante and D. Yao, 145--185,
Springer.

\bibitem{Kobayashi-Miyazawa:2011} Kobayashi, M. and  Miyazawa, M. (2014)
Tail asymptotics of the stationary distribution of a two dimensional reflecting random walk with unbounded upward jumps,
\textit{Adv. in Appl. Probab.},  \textbf{46}(2), 365--399.


\bibitem{Kobayashi-Miyazawa-Zhao:10} Kobayashi, M.,  Miyazawa, M. and Zhao, Y.Q. (2010)
Tail asymptotics of the occupation measure for a Markov additive
process with an M/G/1-type background process,
\textit{Stochastic Models}, \textbf{26}, 463--486.

\bibitem{KST:04} Kroese, D.P., Scheinhardt, W.R.W. and
Taylor, P.G. (2004) Spectral properties of the tandem Jackson
network, seen as a quasi-birth-and-death process, \textit{Annals
of Applied Probability},  \textbf{14}(4), 2057--2089.

\bibitem{Kurkova-Raschel:2011} Kurkova, I.A. and Raschel, K. (2011)
Explicit expression for the generating function counting Gessel's walks,
\textit{Advances in Applied Mathematics}, \textbf{47}, 414--433.

\bibitem{Kurkova-Raschel:2013} Kurkova, I. and Raschel, K. (2013)
Passage time from four to two blocks of opinions in the voter model and walks in the quarter plane,
\textit{Queueing Systems}, \textbf{72}, 219--234.

\bibitem{Kurkova-Suhov:03}
Kurkova, I.A. and  Suhov, Y.M. (2003)
Malyshev's theory and JS-queues. Asymptotics of stationary probabilities,
\textit{The Annals of Applied Probability},
\textbf{13}, 1313--1354.

\bibitem{Li-Miyazawa-Zhao:07} Li, L., Miyazawa, M. and Zhao, Y. (2007) Geometric decay in a QBD process with countable background
states with applications to a join-the-shortest-queue model,
\textit{Stochastic Models}, \textbf{23}, 413--438.

\bibitem{Li-Song-Zhao:2014} Li, H., Song, Y. and Zhao, Y.Q. (2014)
Exact tail asymptotics for the absorbing time --- Revisit of the voter model,
preprint.

\bibitem{Li-Tavakoli-Zhao:11} Li, H., Tavakoli, J. and Zhao, Y.Q. (2013)
Analysis of exact tail asymptotics for singular random walks in
the quarter plane, \textit{Queueing Systems}, \textbf{74}, 151--179.

\bibitem{Li-Zhao:05}
Li, H. and Zhao, Y.Q. (2005) A retrial queue with a constant
retrial rate, server break downs and impatient customers,
\textit{Stochastic Models}, \textbf{21}, 531--550.



\bibitem{Li-Zhao:10b} Li, H. and Zhao, Y.Q. (2011)
Tail asymptotics for a generalized two demand queueing model
--- A kernel method, \textit{Queueing Systems}, \textbf{69}, 77--100.

   \bibitem{Li-Zhao:2018} Li, H. and Zhao, Y.Q. (2018)
A kernel method for exact tail asymptotics --- Random walks in the
quarter plane, \textit{Queueing Models and Service Management}, \textbf{1(1)}, 95--129.

\bibitem{Li-Liu-Zhao:2019} Li, W., Liu, Y. and Zhao Y.Q. (2019)
Exact tail asymptotics for fluid models driven by an $M/M/c$ queue,
\textit{Queueing Systems}, \textbf{91}(3--4), 319--346.

\bibitem{LM:2008} Lieshout, P. and Mandjes, M. (2008)
Asymptotic analysis of L\'{e}vy-driven tandem queues,
\textit{Queueing Systems}, \textbf{60}, 203--226.

\bibitem{Liu-Miyazawa-Zhao:08} Liu, L., Miyazawa, M. and Zhao, Y.Q. (2008)
Geometric decay in level-expanding QBD models, \textit{Annals of
Operations Research}, \textbf{160}, 83--98.

\bibitem{Liu-Wang-Zhao:2016}
Liu, Y., Wang, P. and Zhao, Y.Q. (2016)
Exact stationary tail asymptotics for a Markov modulated two-demand model --- In terms of a kernel method,
\textit{Proceedings of the Ninth International Conference on Matrix-Analytic Methods in Stochastic Models}, June 28--30, Budapest, Hungary, 51--58.

\bibitem{Malyshev:72}
Malyshev, V.A. (1972)
An analytical method in the theory of two-dimensional positive random walks, \textit{Siberian Math. Journal}, \textbf{13}, 1314--1329.

\bibitem{Malyshev:73}
Malyshev, V.A. (1973)
Asymptotic behaviour of stationary probabilities for two dimensional positive random walks, \textit{Siberian Math. Journal}, \textbf{14}, 156--169.

\bibitem{M1977}
Markushevich, A.I. (1977)
\textit{Theory of Funcitons of a Complex Variable}, Vol. I.II.III, English ed.,
translated and edited by Richard A. Silverman, Chelsea Publishing Co., New York.

\bibitem{McDonald:99} McDonald, D.R. (1999) Asymptotics of first
passage times for random walk in an orthant, \textit{Annals of
Applied Probability}, \textbf{9}, 110--145.

\bibitem{Miyazawa:04} Miyazawa, M. (2004) The Markov renewal approach to
$M/G/1$ type queues with countably many background states,
\textit{Queueing Systems}, \textbf{46}, 177--196.

\bibitem{Miyazawa:07} Miyazawa, M. (2007) Doubly QBD process and a solution to the tail decay rate problem,
in \textit{Proceedings of the Second Asia-Pacific Symposium on Queueing Theory and
Network Applications}, Kobe, Japan, 33--42.

\bibitem{Miyazawa:08}
Miyazawa, M. (2009) Tail decay rates in double QBD processes and
related reflected random walks, \textit{Math. OR}, \textbf{34}, 547--575.

\bibitem{Miyazawa:09} Miyazawa, M. (2009) Two sided DQBD process and solutions to the tail decay rate
problem and their applications to the generalized join shortest queue,
in \textit{Advances in Queueing Theory and Network Applications}, edited by W. Yue, Y. Takahashi and H. Takaki, 3--33,
Springer, New York.

\bibitem{Miyazawa:11} Miyazawa, M. (2011) Light tail asymptotics in multidimensional reflecting processes
for queueing networks, \textit{TOP}, \textbf{19(2)}, 233--299.

\bibitem{Miyazawa-Rolski:09} Miyazawa, M. and Rolski, T. (2009)
Tail asymptotics for a L\'{e}vey-driven tandem queue with an
intermediate input, \textit{Queueing Systems}, \textbf{63},
323--353.

\bibitem{Miyazawa-Zhao:04}
Miyazawa, M. and Zhao, Y.Q. (2004) The stationary tail asymptotics
in the $GI/G/1$ type queue with countably many background states,
\textit{Adv. in Appl. Probab.},  \textbf{36}(4), 1231--1251.


\bibitem{Motyer-Taylor:06} Motyer, Allan J. and Taylor, Peter G.
(2006) Decay rates for quasi-birth-and-death process with
countably many phases and tri-diagonal block generators,
\textit{Advances in Applied Probability}, \textbf{38}, 522--544.

\bibitem{Ozawa:2013} Ozawa, T. (2013) Asymptotics for the stationary distribution in a discrete-time two-dimensional quasi-birth-and-death process. \textit{Queueing Systems}, \textbf{74}, 109--149.

\bibitem{Ozawa-Kobayashi:2018}
Ozawa, T. and Kobayashi, M. (2018)
Exact asymptotic formulae of the stationary distribution of a discrete-time two-dimensional QBD process,
\textit{Queueing Systems}, \textbf{90}, 351--403.

\bibitem{Raschel:10} Raschel, K. (2010)
Green functions and Martin compactification for killed random
walks related to SU(3), \textit{Elect. Comm. in Probab.},
\textbf{15}, 176--190.

\bibitem{Song-Liu-Dai:2015} Song, Y., Liu, Z. and Dai, H. (2015) Exact tail asymptotics for a discrete-time preemptive priority queue,
\textit{Acta Mathematicae Applicatae Sinica, English Series}, \textbf{31(1)}, 43--58.

\bibitem{Song-Liu-Zhao:2015} Song, Y., Liu, Z. and Zhao, Y.Q. (2016)
Exact tail asymptotics --- Revisit of a retrial queue with two input streams and two
orbits, \textit{Ann. of Operations Research}, \textbf{247(1)}, 97--120.

\bibitem{Song-Lu:2021} Song, Y. and Lu, H. (2021)
Exact tail asymptotics for the Israeli queue with retrials and non-persistent customers, under review. 

\bibitem{TFM:01} Takahashi, Y., Fujimoto, K. and Makimoto, N. (2001)
Geometric decay of the steady-state probabilities in a
quasi-birth-and-death process with a countable number of phases,
\textit{Stochastic Models}, \textbf{17}(1), 1--24.

\bibitem{Tang-Zhao:08} Tang, J. and Zhao, Y.Q. (2008)
Stationary tail asymptotics of a tandem queue with feedback,
\textit{Annals of Operations Research}, \textbf{160}, 173--189.

\bibitem{W1995}
Williams, R.J. (1995)
Semimartingale reflecting Brownian motions in the orthant,
in \textit{Stochastic Networks. IMA. Vol. Math. Appl.}, \textbf{71}, 125--137, Springer, New York.

\bibitem{W1996}
Williams, R.J. (1996)
On the approximation of queueing networks in heavy traffic,
in \textit{Stochastic Networks: Theory and Applications}, ed. by Kelly, F.P.
Zachary, S. and Ziedins, I., Oxford University Press.

\bibitem{Wright:92} Wright, P. (1992) Two parallel processors with coupled inputs, \textit{Adv. Appl. Prob.}, \textbf{24}, 986--1007.

\bibitem{Ye:2012} Ye, W. (2012) \textit{Longer-Queue-Serve-First System}, M.Sc. Thesis, School of Mathematics and Statistics,
Carleton University. (https://curve.carleton.ca/system/files/theses/29025.pdf)

\bibitem{Ye-Li-Zhao:2015} Ye, W., Li, H. and Zhao, Y.Q. (2015)
Tail behaviour for longest-queue-served-first queueing system ---
A random walk in the half plane, \textit{Stochastic Models}, \textbf{31(3)}, 452--493.

\bibitem{Zafari:2012} Zafari, Z. (2012)
\textit{The Examt Tail Asymptotics Behaviour of the Joint Stationary Distributions of the Generalized-JSQ Model}, master's thesis, University of British Columbia.

\end{thebibliography}
\end{document}